\documentclass[11pt,onecolumn,english]{IEEEtran}
\usepackage[T3,T1]{fontenc}
\usepackage[utf8]{inputenc}
\usepackage{xcolor}
\usepackage{babel}
\usepackage{mathtools}
\usepackage{amsmath}
\usepackage{amsthm}
\usepackage{amssymb}
\usepackage[pdfusetitle,
 bookmarks=true,bookmarksnumbered=true,bookmarksopen=true,bookmarksopenlevel=3,
 breaklinks=false,pdfborder={0 0 1},backref=false,colorlinks=true]
 {hyperref}
\hypersetup{
 pdfborderstyle=,pdfborderstyle={},pdfborderstyle={},pdfborderstyle={},pdfborderstyle={},pdfborderstyle={},pdfborderstyle={},pdfborderstyle={},pdfborderstyle={},pdfborderstyle={},pdfborderstyle={},pdfborderstyle={},pdfborderstyle={},pdfborderstyle={},pdfborderstyle={},pdfborderstyle={},pdfborderstyle={},pdfborderstyle={},pdfborderstyle={},pdfborderstyle={},pdfpagelayout=OneColumn,pdfnewwindow=true,pdfstartview=XYZ,plainpages=false,linkcolor=blue,urlcolor=blue,citecolor=red,anchorcolor=blue,linkcolor=blue,urlcolor=blue,citecolor=red,anchorcolor=blue}

\makeatletter
\usepackage{babel}
\usepackage[mathscr]{eucal}
\usepackage{epsfig,epsf,psfrag}
\usepackage{amssymb,amsmath,amsthm,latexsym}
\usepackage{amsmath,graphicx,xcolor,url}
\usepackage[caption=false]{subfig} 
\usepackage{fixltx2e}%ordering of single and double column floats
\usepackage{array}%array and tabular environments
\usepackage{verbatim}
\usepackage{bm}
\usepackage{algorithmic}
\usepackage{algorithm}
\usepackage{verbatim}
\usepackage{textcomp}
\usepackage{mathrsfs,overpic}
\usepackage{epstopdf}
\usepackage{amsfonts}
\usepackage[numbers]{natbib}
\mathtoolsset{showonlyrefs} 

\usepackage{tikz}
\usepackage{float}
\usepackage{tabularx}
\usepackage{multirow}
\usetikzlibrary{patterns}

\def\1{\mathbf{1}}

\catcode`~=11 \def\UrlSpecials{\do\~{\kern -.15em\lower .7ex\hbox{~}\kern .04em}} \catcode`~=13 

\allowdisplaybreaks[1]

\DeclareMathAlphabet{\mathbsf}{OT1}{cmss}{bx}{n}
\DeclareMathAlphabet{\mathssf}{OT1}{cmss}{m}{sl}% slanted sans serif

\DeclareSymbolFont{bsfletters}{OT1}{cmss}{bx}{n}  
\DeclareSymbolFont{ssfletters}{OT1}{cmss}{m}{n}
\DeclareMathSymbol{\bsfGamma}{0}{bsfletters}{'000}
\DeclareMathSymbol{\ssfGamma}{0}{ssfletters}{'000}
\DeclareMathSymbol{\bsfDelta}{0}{bsfletters}{'001}
\DeclareMathSymbol{\ssfDelta}{0}{ssfletters}{'001}
\DeclareMathSymbol{\bsfTheta}{0}{bsfletters}{'002}
\DeclareMathSymbol{\ssfTheta}{0}{ssfletters}{'002}
\DeclareMathSymbol{\bsfLambda}{0}{bsfletters}{'003}
\DeclareMathSymbol{\ssfLambda}{0}{ssfletters}{'003}
\DeclareMathSymbol{\bsfXi}{0}{bsfletters}{'004}
\DeclareMathSymbol{\ssfXi}{0}{ssfletters}{'004}
\DeclareMathSymbol{\bsfPi}{0}{bsfletters}{'005}
\DeclareMathSymbol{\ssfPi}{0}{ssfletters}{'005}
\DeclareMathSymbol{\bsfSigma}{0}{bsfletters}{'006}
\DeclareMathSymbol{\ssfSigma}{0}{ssfletters}{'006}
\DeclareMathSymbol{\bsfUpsilon}{0}{bsfletters}{'007}
\DeclareMathSymbol{\ssfUpsilon}{0}{ssfletters}{'007}
\DeclareMathSymbol{\bsfPhi}{0}{bsfletters}{'010}
\DeclareMathSymbol{\ssfPhi}{0}{ssfletters}{'010}
\DeclareMathSymbol{\bsfPsi}{0}{bsfletters}{'011}
\DeclareMathSymbol{\ssfPsi}{0}{ssfletters}{'011}
\DeclareMathSymbol{\bsfOmega}{0}{bsfletters}{'012}
\DeclareMathSymbol{\ssfOmega}{0}{ssfletters}{'012}

\newcommand{\qednew}{\nobreak \ifvmode \relax \else
      \ifdim\lastskip<1.5em \hskip-\lastskip
      \hskip1.5em plus0em minus0.5em \fi \nobreak
      \vrule height0.75em width0.5em depth0.25em\fi}

\usepackage{bm,bbm}

\providecommand{\corollaryname}{Corollary}

\providecommand{\lemmaname}{Lemma}
\providecommand{\propositionname}{Proposition}
\providecommand{\remarkname}{Remark}
\providecommand{\theoremname}{Theorem}

\makeatother

\theoremstyle{plain}
\newtheorem{thm}{\protect\theoremname}
\newtheorem{lem}{\protect\lemmaname}
\newtheorem{prop}{\protect\propositionname}
\newtheorem{cor}{\protect\corollaryname}
\theoremstyle{remark}
\newtheorem{rem}{\protect\remarkname}
\providecommand{\corollaryname}{Corollary}
\providecommand{\lemmaname}{Lemma}
\providecommand{\propositionname}{Proposition}
\providecommand{\remarkname}{Remark}
\providecommand{\theoremname}{Theorem}

\begin{document}
\title{On the Concavity of Rényi Entropy and Rényi Entropy Power along
Heat Flows on Manifolds}
\author{Hao Wu and Lei Yu\thanks{H. Wu is with the School of Mathematics and Statistics, Xi'an Jiaotong University, Xi'an 710049, China 
(e-mail: haowu@mail.nankai.edu.cn).
L. Yu is with the School of Statistics and Data Science,
LPMC, KLMDASR, and LEBPS, Nankai University, Tianjin 300071, China
(e-mail: leiyu@nankai.edu.cn). This
work was supported by the National Key Research and Development Program
of China under grant 2023YFA1009604 and the NSFC under grant 62101286. The work of H. Wu was also supported by the China Postdoctoral Science Foundation under grant 2026M793399.
%(\textit{The corresponding author: Lei Yu.})
}}
\maketitle
\begin{abstract}
In this paper, we establish the concavity of Rényi entropy and Rényi
entropy power along heat flow on Euclidean space of arbitrary dimension,
for Rényi orders $\alpha\in[1,2]$. More generally,
we prove the corresponding result along the heat flow generated by
the Witten Laplacian on complete Riemannian manifolds under a suitable
curvature-dimension condition.  Our argument builds on novel functional
inequalities that were recently developed for the study of Tsallis
entropy concavity along the heat flow. We extend both these inequalities
and the resulting concavity theorem from the Euclidean setting to
the Riemannian manifold setting. In addition, we establish the complete
monotonicity of Tsallis entropy of order $2$ along the heat flow
and derive its spectral decomposition.
\end{abstract}


\section{Introduction}

The evolution of entropy along heat flow provides a fundamental link
between information theory, probability, and geometric analysis.
For Shannon entropy, this connection appears through the de Bruijn
identity, Fisher information, and the concavity of entropy power under
Gaussian smoothing \cite{McKean1966,costa1985new,dembo1989simple,villani2000short}.
It is natural to ask which of these properties remain valid when
Shannon entropy is replaced by R\'enyi or Tsallis entropy and when
the Euclidean Laplacian is replaced by a diffusion generator on a
weighted Riemannian manifold.

\subsection{Geometric setting and entropy functionals}

Let $(M,g)$ be an $n$-dimensional complete Riemannian manifold,
$Ric$ the Ricci curvature of $g$, $\phi\in C^{2}(M)$ and $d\mu=e^{-\phi}d\omega$,
where $\omega$ is the Riemannian volume measure on $(M,g)$. The
Witten Laplacian acting on smooth functions is defined by
\[
L=\Delta-\nabla\phi\cdot\nabla.
\]
Denote the set of all infinitely differentiable functions on $M$
whose support set is compact as $C^{\infty}_{0}(M)$, and the set
of all $k$-times continuously differentiable, bounded functions defined
on $M$ as $C^{k}_{b}(M)$. For any $u,v\in C^{\infty}_{0}(M)$, the
integration-by-parts formula holds
\[
\int_{M}\langle\nabla u,\nabla v\rangle d\mu=-\int_{M}Luvd\mu=-\int_{M}uLvd\mu.
\]

The diffusion product rule
\[
L(uv)=uLv+vLu+2\langle\nabla u,\nabla v\rangle
\]
and the corresponding chain rule are standard; see
\cite{bakry2006diffusions,bakry2013analysis}.

The generalized Bochner--Weitzenböck formula holds \cite{bakry2006diffusions}
\[
\frac{1}{2}L|\nabla u|^{2}=\langle\nabla u,\nabla Lu\rangle+|\nabla^{2}u|^{2}+\langle Ric(L)\nabla u,\nabla u\rangle,
\]
where
\[
Ric(L):=Ric+\nabla^{2}\phi,
\]
and $|\nabla^{2}u|$ denotes the Hilbert-Schmidt norm of the Hessian
of $u$, i.e., $\nabla^{2}u$. Let $m\ge n$ be a constant. The $m$-dimensional
Bakry--Emery Ricci curvature associated with the Witten Laplacian
$L$ is defined by \cite{bakry2006diffusions,lott2002some,li2005liouville}
\[
Ric_{m,n}(L):=Ric+\nabla^{2}\phi-\frac{\nabla\phi\otimes\nabla\phi}{m-n}.
\]
Here we adopt the convention that $m=n$ if and only if $\phi$ is
a constant. 

The heat kernel $p_{t}(x,y)$ of the Witten Laplacian $L$ is the
fundamental solution to the heat equation
\[
\partial_{t}u=Lu.
\]
We shall also consider general positive solutions $u=u(x,t)$ to the
same equation. Throughout this paper, both fundamental solutions and
positive solutions are assumed to be normalized with respect to $\mu$,
namely,
\[
\int_{M}u(x,t)d\mu(x)=1,\quad\forall t>0.
\]

For $\alpha\in(0,+\infty)\backslash\{1\}$, the Tsallis entropy of
order $\alpha$ for a probability distribution with density $u$ is
defined by
\[
S_{\alpha}(u):=\frac{1}{1-\alpha}\left(\int_{M}u^{\alpha}(x)d\mu(x)-1\right),
\]
and the Rényi entropy of order $\alpha$ is defined by
\[
h_{\alpha}(u):=\frac{1}{1-\alpha}\log\left(1+(1-\alpha)S_{\alpha}(u)\right)=\frac{1}{1-\alpha}\log\int_{M}u^{\alpha}(x)d\mu(x).
\]
For $\alpha=1$, the Tsallis and Rényi entropies are defined by continuous
extension. Specifically, for $\alpha=1$, the Tsallis and Rényi entropy
 both reduce to the Shannon entropy
\[
h(u)=-\int_{M}u(x)\log u(x)d\mu(x).
\]

\subsection{Heat-flow inequalities and related literature}

In \cite{li2012perelman}, Li computed the first and second time derivatives
of Shannon entropy  along the heat flow generated by the Witten Laplacian
on complete Riemannian manifolds, thereby extending the corresponding
results for the standard Laplacian on Euclidean spaces. The Euclidean
case is more closely connected with information theory and with several
important topics, including the entropy power inequality (EPI) \cite{Shannon48,StamEPI59}.
In particular, the de Bruijn identity states that, for a positive
solution of the standard heat equation,
\[
\partial_{t}u=\Delta u
\]
on $\mathbb{R}^{n}$, we have
\[
\frac{d}{dt}h(u)=I(u),
\]
where $I(u)=\int|\nabla u|^{2}/udx$ is the Fisher information \cite{fisher1925theory,StamEPI59}.

The de Bruijn identity plays an important role in the proof of EPI
and other information-theoretic inequalities. Further investigation
of the second time derivative of entropy along the heat flow has led
to refined inequalities, including Costa’s EPI \cite{costa1985new},
which strengthens the classical EPI when one of the two random variables
involved in the EPI is Gaussian. 

R\'enyi entropy, introduced in \cite{Renyi1961OnMO}, is also closely
related to the $L^{\alpha}$-norm of a probability density. This relationship
connects its evolution to contractivity and hypercontractivity of Markov
semigroups, logarithmic Sobolev inequalities, and the Bakry--\'Emery
approach to diffusion inequalities
\cite{gross1975logarithmic,bakry2006diffusions,bakry2013analysis}.
Moreover, the order-two quantity $\int_M(P_t f)^2d\mu$ is the noise
stability of $f$ with respect to the semigroup, a functional closely
connected to isoperimetric inequalities and other geometric estimates
\cite{ledoux1994semigroup,wu2025onthecomp}.

Beyond the first two entropy derivatives, McKean's questions on the
behavior of entropy along heat flow led to investigations of higher
order alternating-sign inequalities. Cheng and Geng proved results on
the third- and fourth-order derivatives of Shannon entropy
\cite{Cheng2015Higher}, and subsequent work has studied related
properties of Fisher information
\cite{ledoux2021log,ledoux2022differentials}.
Hung investigated the corresponding problem for Tsallis entropy
\cite{hung2022generalization}, whereas Wu, Yu, and Guo examined the
R\'enyi case, its entropy power, and order-two Tsallis entropy
\cite{wu2025onthecomp}. The distinction between the two entropies is
important: the logarithm in R\'enyi entropy introduces products of
integrals into higher derivatives, so inequalities for Tsallis entropy
cannot automatically be transferred to R\'enyi entropy.

A different and extensively studied direction concerns the concavity of
R\'enyi entropy power along \emph{nonlinear} diffusion equations. It
includes the work of Savar\'e and Toscani on Euclidean spaces and its
extensions by Guo, Yuan, and Gao and by Li and Li to Riemannian
manifolds \cite{savare2014concavity,guo2021generalization,li2020renyi}.
The present paper instead considers the \emph{linear} Witten heat
semigroup, where the R\'enyi order varies but the diffusion generator
does not. Our aim is to understand how curvature assumptions enter the
second-order differential inequalities in this setting.

\subsection{Contributions and the role of curvature}

A recent paper \cite{sun2026concavity} studied the second time derivative
of Tsallis entropy along the heat flow on $n$-dimensional Euclidean
space and proved that it is nonpositive when the Tsallis order $1\le\alpha\le3$,
thereby extending the corresponding one-dimensional results of \cite{hung2022generalization,wu2025onthecomp}.
In the present paper, we investigate the first and second time derivatives
of Tsallis and Rényi entropy along the heat flow generated by the
Witten Laplacian on complete Riemannian manifolds. The main contributions
of this paper are twofold. First, we extend the concavity of Tsallis
entropy in \cite{sun2026concavity} from the standard heat equation
on Euclidean spaces to the heat equation associated with the Witten
Laplacian on complete, possibly noncompact, Riemannian manifolds satisfying
suitable curvature-dimension conditions. Moreover, in this geometric
setting, we establish the concavity of the Rényi entropy and a suitably
defined Rényi entropy power along the Witten heat flow for Rényi orders
$1\le\alpha\le2$. In particular, this provides a higher-dimensional counterpart
to the corresponding one-dimensional result in \cite{wu2025onthecomp}.

%\subsubsection{Curvature assumptions and the order-two endpoint}

A basic question is which conclusions depend on the geometry of the
underlying manifold. The weighted Bochner identity contains the
Bakry--\'Emery curvature tensor $Ric(L)$, whereas a finite-dimensional
entropy-power estimate additionally exploits the curvature--dimension
condition $Ric_{m,n}(L)\ge0$. In this paper these conditions enter the
proofs through explicit integration-by-parts identities and functional
inequalities. In contrast, the alternating derivatives of the order-two
Tsallis entropy follow from the spectral representation of a symmetric
Markov semigroup and do not intrinsically require a curvature bound.
This separation provides another motivation for considering the
weighted-manifold setting.

\subsection{Notation}

Throughout this paper, for any $\beta>0$, we write $\Vert f\Vert_{\beta}=\left(\int_{M}|f(x)|^{\beta}d\mu\right)^{\frac{1}{\beta}}$,
and $\Vert f\Vert_{\infty}=\text{esssup}_{x\in M}|f(x)|$. We denote
by $L^{\beta}(M,\mu)$ the space of all measurable $f$ such that
$\Vert f\Vert_{\beta}<\infty$. When $0<\beta<1$, $\Vert\cdot\Vert_{\beta}$
is understood as a quasi-norm rather than a norm. Here $f$ may be
a scalar function, a vector field, or, more generally, a tensor field.
Accordingly, $|f|$ denotes the pointwise absolute value in the scalar
case and the pointwise norm induced by the Riemannian metric in the
vector- or tensor-valued case. With this notation, $S_{\alpha}(u)=\frac{1}{1-\alpha}\left(\Vert u\Vert^{\alpha}_{\alpha}-1\right)$,
$h_{\alpha}(u)=\frac{\alpha}{1-\alpha}\log\Vert u\Vert_{\alpha}$.
For a Riemannian manifold $(M,g)$, $\text{Riem}$ denotes the Riemann
curvature tensor, and $\nabla^{k}\text{Riem}$ denotes its $k$-th
covariant derivative with respect to the Levi-Civita connection. 

\subsection{Main results}

The key theorems are collected here; their proofs, together with the auxiliary functional inequalities, are given in the subsequent sections. We keep the hypotheses for fundamental solutions and general positive solutions separate throughout.

\subsubsection{Dissipation identities and regularity assumptions}

The following dissipation formulas record the precise regularity hypotheses used in our manifold statements. The heat-kernel version applies to the fundamental solution; the positive-solution version uses the local-in-time integrability assumption \eqref{eq:2-9-2-4}.

\begin{thm}
\label{thm:derivatives of TsallisEnt-1}Let $(M,g)$ be an $n$-dimensional
complete Riemannian manifold. Assume that $\Vert\nabla^{j}Riem\Vert_{\infty}<C^{(1)}_{j}$
for some constant $C^{(1)}_{j}>0,0\le j\le3$. Let $\phi\in C^{4}(M)$
be such that $\nabla\phi\in C^{3}_{b}(M)$, and suppose that there
exist constants $m\ge n$ and $K\ge0$ such that
\[
Ric_{m,n}(L)\ge-K.
\]
Then, for the fundamental solution $u(x,t)$ to the heat equation
\[
\partial_{t}u=Lu,
\]
and any $\alpha\ge1$, we have
\begin{equation}
\frac{d}{dt}S_{\alpha}(u)=\alpha\int_{M}|\nabla u(t)|^{2}u^{\alpha-2}(t)d\mu,\label{eq:3-12-15}
\end{equation}
and
\begin{equation}
\frac{d^{2}}{dt^{2}}S_{\alpha}(u)=-\alpha\left(2\int_{M}u^{\alpha-2}(Lu)^{2}d\mu+(\alpha-2)\int_{M}u^{\alpha-3}|\nabla u|^{2}Lud\mu\right)\label{eq:2-9-2}
\end{equation}
for every $t>0$.
\end{thm}

\begin{thm}
\label{thm:Let--be}Fix an order $\alpha\ge1$. Let $(M,g)$ be an $n$-dimensional complete Riemannian
manifold, and $\phi\in C^{2}(M)$. Let $u$ be a smooth positive solution
of the heat equation
\[
\partial_{t}u=Lu.
\]
Assume that, for every $0<a<b<\infty$,
\begin{equation}
\sup_{a\le t\le b}\int_{M}u^{\alpha}d\mu<\infty,\label{eq:2-9-2-3}
\end{equation}
and
\begin{equation}
\int^{b}_{a}\int_{M}u^{\alpha-2}\left(|Lu|^{2}+|\nabla Lu|^{2}+\frac{|\nabla u|^{4}}{u^{2}}\right)d\mu\,dt<\infty.\label{eq:2-9-2-4}
\end{equation}
At the Shannon endpoint $\alpha=1$, assume in addition that
\[
\sup_{a\le t\le b}\int_Mu(t)(\log u(t))^2\,d\mu<\infty
\]
for every $0<a<b<\infty$, so that the entropy and the logarithmic
cutoff terms are finite.
Then, for the fixed order $\alpha\ge1$,
\begin{equation}
\frac{d}{dt}S_{\alpha}(u)=\alpha\int_{M}|\nabla u(t)|^{2}u^{\alpha-2}(t)d\mu\label{eq:3-12-15-1}
\end{equation}
and
\begin{equation}
\frac{d^{2}}{dt^{2}}S_{\alpha}(u)=-\alpha\left(2\int_{M}u^{\alpha-2}(Lu)^{2}d\mu+(\alpha-2)\int_{M}u^{\alpha-3}|\nabla u|^{2}Lud\mu\right)\label{eq:2-9-2-2}
\end{equation}
hold for almost every $t>0$. Furthermore, if the integrals in \eqref{eq:3-12-15-1}
and \eqref{eq:2-9-2-2} are convergent and continuous for every $t>0$,
then \eqref{eq:3-12-15-1} and \eqref{eq:2-9-2-2} hold for every
$t>0$.
\end{thm}

\subsubsection{Concavity of Rényi entropy}

We state the Rényi concavity theorem first for the fundamental solution, then for positive solutions and Euclidean heat flow.

\begin{thm}
\label{thm:Under-the-assumptions-1}Under the assumptions of Theorem
\ref{thm:derivatives of TsallisEnt-1}, furthermore, if
\[
Ric(L)\ge0,
\]
then, for $\alpha\in[1,2]$, 
\[
\frac{d^{2}}{dt^{2}}h_{\alpha}(u_{t})\le0.
\]
\end{thm}

\begin{thm}
\label{thm:Renyi-positive-concavity}
Under the assumptions of Theorem \ref{thm:Let--be}, suppose further
that the integrals in \eqref{eq:3-12-15-1} and \eqref{eq:2-9-2-2}
are convergent and continuous for every $t>0$, that
\[
Ric(L)\ge0,
\]
and that the following integration by parts formulas
 \begin{align}
\int_{M}wLwd\mu & =-\int_{M}|\nabla w|^{2}d\mu,\label{eq:4-1-1-1}\\
\int_{M}L|\nabla w|^{2}d\mu & =0,\label{eq:4-1-2-1}\\
\int_{M}\langle\nabla Lw,\nabla w\rangle d\mu & =-\int_{M}|Lw|^{2}d\mu,\label{eq:4-1-3-1}\\
\int_{M}Lw|\nabla\sqrt{w}|^{2}d\mu & =-\int_{M}\langle\nabla w,\nabla|\nabla\sqrt{w}|^{2}\rangle d\mu\label{eq:4-1-4-1}
\end{align} hold with $w=u_t^{\alpha/2}$ for every $\alpha\in[1,2]$ and $t>0$. Then, for every $\alpha\in[1,2]$ and $t>0$, 
\[
\frac{d^{2}}{dt^{2}}h_{\alpha}(u_{t})\le0.
\]
\end{thm}

\begin{thm}
\label{thm:For-any-positive}For any positive integer $n$, let $u_{t}$
be a positive normalized solution of the heat equation
\[
\partial_{t}u=\Delta u
\]
on $n$-dimensional Euclidean space obtained by Gaussian convolution
from nonnegative $L^{1}$ initial data. Then
\[
\frac{d^{2}}{dt^{2}}h_{\alpha}(u_{t})\le0
\]
for $\alpha\in[1,2]$.
\end{thm}

\subsubsection{Concavity of Rényi entropy power}

The finite-dimensional curvature condition also yields a quantitative entropy-power inequality. The effective dimension $m$ is essential in the weighted setting.

\begin{thm}
\label{thm:Under-the-assumptions-3}Under the assumptions of Theorem
\ref{thm:derivatives of TsallisEnt-1}, furthermore, if  $\alpha\in[1,2]$,
$0<\beta\le2(2-\alpha)$,   and
\[
Ric_{m,n}(L)\ge0,
\]
we have
\[
\frac{d^{2}}{dt^{2}}e^{\frac{\beta}{m}h_{\alpha}(u)}\le0.
\]
\end{thm}

\begin{thm}
\label{thm:Renyi-positive-power}
Under the assumptions of Theorem \ref{thm:Let--be}, suppose further
that the integrals in \eqref{eq:3-12-15-1} and \eqref{eq:2-9-2-2}
are convergent and continuous for every $t>0$, that
\[
Ric_{m,n}(L)\ge0\quad\text{for some finite }m\ge n,
\]
and that the integration by parts formulas \eqref{eq:4-1-1-1}, \eqref{eq:4-1-2-1}, \eqref{eq:4-1-3-1}, and \eqref{eq:4-1-4-1} in Theorem \ref{thm:Renyi-positive-concavity}
hold with $w=u_t^{\alpha/2}$ for every $\alpha\in[1,2]$ and $t>0$. Then, for all $\alpha\in[1,2]$, $t>0$, and $0<\beta\le2(2-\alpha)$,
\[
\frac{d^{2}}{dt^{2}}e^{\frac{\beta}{m}h_{\alpha}(u)}\le0.
\]
\end{thm}

\begin{thm}
\label{thm:Renyi-euclidean-power}
For any positive integer $n$, let $u_t$ be a positive normalized solution
of the heat equation
\[
\partial_{t}u=\Delta u
\]
on $n$-dimensional Euclidean space obtained by Gaussian convolution
from nonnegative $L^{1}$ initial data. Then, for $\alpha\in[1,2]$,
$0<\beta\le2(2-\alpha)$ and $t>0$,
\[
\frac{d^{2}}{dt^{2}}e^{\frac{\beta}{n}h_{\alpha}(u_t)}\le0.
\]
\end{thm}

\subsubsection{Concavity of Tsallis entropy}

For comparison, the same geometric approach extends the range $1\le\alpha\le3$ for Tsallis entropy to the Witten heat flow.

\begin{thm}
\label{thm:Under-the-assumptions}Under the assumptions of Theorem
\ref{thm:derivatives of TsallisEnt-1}, furthermore, if
\[
Ric(L)\ge0,
\]
then, for $\alpha\in[1,3]$ and $t>0$, 
\[
\frac{d^{2}}{dt^{2}}S_{\alpha}(u_{t})\le0.
\]
\end{thm}

\begin{thm}
\label{thm:Under-the-assumptions-2}Under the assumptions of Theorem
\ref{thm:Let--be}, suppose further that the integrals in \eqref{eq:3-12-15-1} and \eqref{eq:2-9-2-2}
are convergent and continuous for every $t>0$, that
\[
Ric(L)\ge0,
\]
and that the integration by parts formulas \eqref{eq:4-1-2-1}, \eqref{eq:4-1-3-1}, and \eqref{eq:4-1-4-1} in Theorem \ref{thm:Renyi-positive-concavity}
hold with $w=u_t^{\alpha/2}$ for every $\alpha\in[1,3]$
and $t>0$. Then, for every $\alpha\in[1,3]$
and $t>0$,
\[
\frac{d^{2}}{dt^{2}}S_{\alpha}(u_{t})\le0.
\]
\end{thm}

\subsubsection{Complete monotonicity at order two}

At the order-two endpoint the argument is spectral and therefore fundamentally different: curvature is not intrinsically needed for complete monotonicity of the $L^2$-energy of a symmetric Markov semigroup.

\begin{thm}
\label{thm:fundamental-CM}
Let $(M,g)$ be an $n$-dimensional complete Riemannian manifold.
 Assume that $\Vert\nabla^{j}Riem\Vert_{\infty}<C^{(1)}_{j}$ for
positive constants $C^{(1)}_{j}>0$ for every integer $j\ge0$. Let $\phi\in C^{\infty}(M)$
be such that $\nabla\phi\in C^{\infty}_{b}(M)$, and suppose that
there exist two constants $m\ge n$ and $K\ge0$ such that
\[
Ric_{m,n}(L)\ge-K.
\]

Then, for the fundamental solution $u(x,t)$ to the heat equation
\[
\partial_{t}u=Lu,
\]
we have that $\frac{d}{dt}S_{2}(u(t))$ is completely monotone in $t\in(0,+\infty)$.
\end{thm}

\begin{thm}
\label{thm:positive-CM}
Let $(M,g)$ be an $n$-dimensional complete
Riemannian manifold, and $\phi\in C^{2}(M)$. Let $u$ be a smooth
positive solution of the heat equation
\[
\partial_{t}u=Lu.
\]
Assume that, for every $0<a<b<\infty$,
\begin{equation}
\sup_{t\in[a,b]}\int_{M}u^{2}d\mu<\infty,
\end{equation}
and
\begin{equation}
\int^{b}_{a}\int_{M}\left(|L^{k}u|^{2}+|\nabla L^{k}u|^{2}\right)d\mu dt<\infty
\end{equation}
for every nonnegative integer $k$. Assume further that, for every nonnegative integer $k$, the integrals
\[
\int_{M}|\nabla L^{k}u|^{2}d\mu
\quad\text{and}\quad
\int_{M}(L^{k+1}u)^{2}d\mu
\]
converge and are
continuous with respect to $t$ for every $t>0$. Then we have that $\frac{d}{dt}S_{2}(u(t))$ is completely monotone
in $t\in(0,+\infty)$.
\end{thm}

\begin{thm}
\label{thm:For-any-positive-1}For any positive integer $n$, let
$u(x,t)$ be a positive normalized solution of the heat equation
\[
\partial_{t}u=\Delta u
\]
on $n$-dimensional Euclidean space obtained by Gaussian convolution
from nonnegative $L^{1}$ initial data, we have that $\frac{d}{dt}S_{2}(u(t))$ is
completely monotone in $t\in(0,+\infty)$.
\end{thm}

\subsection{Organization of the proofs}

The remainder of this paper is organized around the proofs of the results
stated above. Section \ref{sec:Tsallis-entropy-dissipation} establishes
the dissipation identities under the stated heat-kernel and positive-solution
assumptions. Section \ref{sec:Preliminaries-for-the} derives the auxiliary
integration-by-parts identities and curvature-dependent functional
inequalities used in the concavity arguments. Section
\ref{sec:The-concavity-of} proves the Tsallis concavity theorems. Section
\ref{sec:Tsallis-entropy-of} proves the order-two complete-monotonicity
results and gives the spectral interpretation. Sections
\ref{sec:R=0000E9nyi-entropy-dissipation} and
\ref{sec:Concavity-of-R=0000E9nyi} prove, respectively, the R\'enyi
entropy and R\'enyi entropy-power concavity theorems. The appendix
collects the heat-kernel estimates required for the noncompact case.

\section{\label{sec:Tsallis-entropy-dissipation}Proofs of the Tsallis entropy dissipation theorems}

We compute the first and second time derivatives of Tsallis entropy
$S_{\alpha}$ along the heat flow, for orders $\alpha\ge1$.

\subsection{Proof for the fundamental solution}

\begin{IEEEproof}
When $\alpha=1,$ $S_{\alpha}$ reduces to the Shannon entropy. An
equivalent formulation of the result is given in Theorems 5.1 and
5.2 of \cite{li2016hamilton}, so we do not consider this case here.
It remains only to prove the case $\alpha>1$. Let $\{\xi_{k}\}^{\infty}_{k=1}$
be an increasing sequence of functions in $C^{\infty}_{0}(M)$ such
that $0\le\xi_{k}\le1$ and $|\nabla\xi_{k}|\le\frac{1}{k}$. Fix
an arbitrary point $y\in M$ and let $u(x,t)=p_{t}(x,y)$ be the fundamental
solution. Fix any $0<a<b<\infty$, and let $a\le s<t\le b$.

First, we compute $\frac{d}{dt}S_{\alpha}(u)$. Let
\[
H_{k}(t):=\int_{M}u^{\alpha}(t)\xi_{k}d\mu,\quad H(t):=\int_{M}u^{\alpha}(t)d\mu.
\]
 Since $\xi_{k}\in C^{\infty}_{0}(M)$ and is independent of $t$,
taking derivative with respect to $t$, we have
\begin{align}
H_{k}'(t) & =\frac{d}{dt}H_{k}(t)\nonumber \\
 & =\int_{M}\frac{d}{dt}(u^{\alpha}(t))\xi_{k}d\mu\nonumber \\
 & =\alpha\int_{M}u^{\alpha-1}(t)Lu(t)\xi_{k}d\mu.\label{eq:2-1-2}
\end{align}
Applying integrations by parts, we obtain
\begin{align}
H_{k}'(t) & =\alpha\int_{M}u^{\alpha-1}(t)Lu(t)\xi_{k}d\mu\nonumber \\
 & =-\alpha\int_{M}\langle\nabla u(t),\nabla(u^{\alpha-1}(t)\xi_{k})\rangle d\mu\nonumber \\
 & =-\alpha(\alpha-1)\int_{M}|\nabla u(t)|^{2}u^{\alpha-2}(t)\xi_{k}d\mu-\alpha\int_{M}\langle\nabla u(t),\nabla\xi_{k}\rangle u^{\alpha-1}(t)d\mu\nonumber \\
 & =:-\alpha(\alpha-1)I_{1,k}(t)-\alpha I_{2,k}(t).\label{eq:2-2-2}
\end{align}
Thus, 
\begin{equation}
H_{k}(t)-H_{k}(s)=-\alpha(\alpha-1)\int^{t}_{s}I_{1,k}(r)dr-\alpha\int^{t}_{s}I_{2,k}(r)dr.\label{eq:2-12-2}
\end{equation}
By the monotone convergence theorem, it is easy to see that
\begin{equation}
H_{k}(t)\uparrow H(t)\label{eq:2-10-2}
\end{equation}
and
\begin{equation}
I_{1,k}(t)\uparrow I_{1}(t):=\int_{M}|\nabla u(t)|^{2}u^{\alpha-2}(t)d\mu\label{eq:2-13-2}
\end{equation}
as $k\to\infty$. Moreover, by Proposition 3.1 in \cite{li2016hamilton},
for every fixed $t>0$ and $y\in M$, the heat kernel $u(\cdot,t)=p_{t}(\cdot,y)$
satisfies a Gaussian upper bound with at most polynomial and linear-exponential
factors. Hence
\begin{equation}
\sup_{a\le t\le b}\Vert u(\cdot,t)\Vert_{\infty}<+\infty.\label{eq:3-12-1}
\end{equation}
By Lemma \ref{thm:gradient estimates-1}, there exists a constant
$C_{0}(K,m,n,b)>0$, such that
\begin{equation}
|\nabla u(x,t)|\le C_{0}(K,m,n,b)\left(\frac{d(x,y)}{t}+\frac{1}{\sqrt{t}}\right)u(x,t).\label{eq:2-5-3}
\end{equation}
It follows that
\begin{align*}
\int_{M}|\nabla u(t)|^{2}u^{\alpha-2}(t)d\mu & =\int_{M}|\nabla u(t)|^{2}u^{-1}(t)u^{\alpha-1}(t)d\mu\\
 & \le\Vert u(\cdot,t)\Vert^{\alpha-1}_{\infty}\int_{M}\frac{|\nabla u(t)|^{2}}{u}d\mu\\
 & \le C_{0}(K,m,n,b)^{2}\Vert u(\cdot,t)\Vert^{\alpha-1}_{\infty}\int_{M}\left(\frac{d(x,y)}{t}+\frac{1}{\sqrt{t}}\right)^{2}u(t)d\mu.
\end{align*}
By Lemma \ref{lem:Let-.-Let} and \eqref{eq:3-12-1},  we have
\[
\int_{M}|\nabla u(t)|^{2}u^{\alpha-2}(t)d\mu<\infty,\quad\forall t>0,
\]
and
\begin{equation}
\sup_{a\le t\le b}I_{1}(t)<\infty.\label{eq:3-12-3}
\end{equation}
By \eqref{eq:2-13-2} and \eqref{eq:3-12-3}, we have 
\begin{equation}
\int^{t}_{s}I_{1,k}(r)dr\to\int^{t}_{s}I_{1}(r)dr\label{eq:3-12-6}
\end{equation}
as $k\to\infty$, and
\begin{equation}
I_{1}(t)\in L^{1}([a,b]).\label{eq:3-12-9}
\end{equation}

Since $|\nabla\xi_{k}|\le\frac{1}{k}$,
\begin{align}
|I_{2,k}(t)| & =|\int_{M}\langle\nabla u(t),\nabla\xi_{k}\rangle u^{\alpha-1}(t)d\mu|\nonumber \\
 & \le\frac{1}{k}\int_{M}|\nabla u(t)|u^{\alpha-1}(t)d\mu\nonumber \\
 & =\frac{1}{k}\int_{M}|\nabla u(t)|u^{\frac{\alpha-2}{2}}(t)u^{\frac{\alpha}{2}}(t)d\mu\nonumber \\
 & \le\frac{1}{k}\left(\int_{M}|\nabla u(t)|^{2}u^{\alpha-2}(t)d\mu\right)^{1/2}\left(\int_{M}u^{\alpha}(t)d\mu\right)^{1/2},\label{eq:3-12-16}
\end{align}
where the last inequality comes from Hölder's inequality. By the
fact that $\int_{M}u(t)d\mu=1,$ and \eqref{eq:3-12-1},
\begin{align}
\sup_{a\le t\le b}\int_{M}u^{\alpha}(t)d\mu & =\sup_{a\le t\le b}\int_{M}u^{\alpha-1}(t)u(t)d\mu\nonumber \\
 & \le\sup_{a\le t\le b}\Vert u(\cdot,t)\Vert^{\alpha-1}_{\infty}\nonumber \\
 & <\infty.\label{eq:3-12-4}
\end{align}
Then,  by \eqref{eq:3-12-16}\eqref{eq:3-12-3}\eqref{eq:3-12-4}
and dominated convergence theorem,
\begin{equation}
I_{2,k}(t)\to0,\label{eq:3-12-11}
\end{equation}
and
\begin{equation}
\int^{t}_{s}I_{2,k}(r)dr\to0\label{eq:3-12-7}
\end{equation}
as $k\to\infty$. Then taking $k\to\infty$ on both sides of \eqref{eq:2-12-2},
and by \eqref{eq:2-10-2}\eqref{eq:3-12-6}\eqref{eq:3-12-7}, 
\begin{equation}
H(t)-H(s)=-\alpha(\alpha-1)\int^{t}_{s}I_{1}(r)dr.\label{eq:3-12-8}
\end{equation}
Taking derivative with respect to $t$ on both side of \eqref{eq:3-12-8},
and note that \eqref{eq:3-12-9} and $a<b$ are arbitrary in $(0,\infty)$,
we derive that
\begin{equation}
H'(t)=-\alpha(\alpha-1)I_{1}(t)\label{eq:3-12-12}
\end{equation}
holds for almost every $t>0$.

Equivalently, \eqref{eq:3-12-15} holds for almost every $t>0$. To
show that it holds for every $t>0$, it suffices to prove that $I_{1}(t)$
is continuous on $(0,\infty)$. Indeed, set $r=r(x)=d(x,y)$. Since
$M$ is complete, the closed ball $\overline{B(y,R)}=\left\{ x\in M,r(x)\le R\right\} $
is compact. Moreover, the heat kernel $u(x,t)=p_{t}(x,y)$ is smooth
and strict positive on $M\times(0,\infty)$. Consequently, $|\nabla u(x,t)|^{2}u^{\alpha-2}(x,t)$
is continuous on the compact set $\overline{B(y,R)}\times[a,b]$,
and hence is bounded there. Set $\mathcal{I}_{R}(t)=\int_{\left\{ r\le R\right\} }|\nabla u(x,t)|^{2}u^{\alpha-2}(x,t)d\mu(x)$.
It follows from the dominated convergence theorem that
\[
t\mapsto\mathcal{I}_{R}(t)
\]
 is continuous on $[a,b]$.

It remains to show that 
\[
\mathcal{I}_{R}(t)\to I_{1}(t)
\]
uniformly on $[a,b]$ as $R\to\infty$. 

Estimate \eqref{eq:2-5-3} gives
\[
|\nabla u(x,t)|^{2}u^{\alpha-2}(x,t)\le C_{0}(K,m,n,b)^{2}\left(\frac{r}{t}+\frac{1}{\sqrt{t}}\right)^{2}u^{\alpha}(x,t),
\]
since $t\ge a>0$, this implies that there exists a constant $C=C(K,m,n,a,b)$,
such that
\[
|\nabla u(x,t)|^{2}u^{\alpha-2}(x,t)\le C\left(1+r^{2}\right)u^{\alpha}(x,t).
\]
Choose an integer $Q>2$. On the set $\left\{ r>R\right\} $,
\[
1\le R^{-Q}r^{Q}
\]
and
\[
r^{2}\le R^{2-Q}r^{Q}.
\]
Therefore,
\begin{align*}
\left|\mathcal{I}_{R}(t)-I_{1}(t)\right| & =\int_{\left\{ r>R\right\} }|\nabla u(x,t)|^{2}u^{\alpha-2}(x,t)d\mu(x)\\
 & \le C\int_{\left\{ r>R\right\} }\left(1+r^{2}\right)u^{\alpha}(x,t)d\mu(x)\\
 & \le C\sup_{a\le t\le b}\Vert u(\cdot,t)\Vert^{\alpha-1}_{\infty}\int_{\left\{ r>R\right\} }\left(1+r^{2}\right)u(x,t)d\mu(x)\\
 & \le C\left(R^{-Q}+R^{2-Q}\right)\sup_{a\le t\le b}\Vert u(\cdot,t)\Vert^{\alpha-1}_{\infty}\int_{M}r(x)^{Q}u(x,t)d\mu(x).
\end{align*}
By \eqref{eq:3-12-1} and Lemma \ref{lem:Let-.-Let},   and since
$Q>2$, it follows that
\[
\lim_{R\to\infty}\sup_{a\le t\le b}\left|\mathcal{I}_{R}(t)-I_{1}(t)\right|=0.
\]
Thus $\mathcal{I}_{R}(t)\to I_{1}(t)$ uniformly on $[a,b]$. Since
each $\mathcal{I}_{R}(t)$ is continuous on $[a,b]$, the uniform
limit $I_{1}(t)$ is also continuous on $[a,b]$. As the interval
$[a,b]\subset(0,\infty)$ was arbitrary, $I_{1}(t)$ is continuous
in $t>0$.

 Next, we compute $\frac{d^{2}}{dt^{2}}S_{\alpha}(u)$. By \eqref{eq:2-2-2}\eqref{eq:2-13-2}
and \eqref{eq:3-12-11},
\begin{equation}
H'_{k}(t)\to H'(t)\label{eq:2-5-11}
\end{equation}
as $k\to\infty$. Taking derivative with respect to $t$ on both sides
of \eqref{eq:2-1-2}, and applying integration by parts, we have
\begin{align*}
H''_{k}(t) & =\frac{d}{dt}H'_{k}(t)\\
 & =\alpha\frac{d}{dt}\int_{M}u^{\alpha-1}(t)Lu(t)\xi_{k}d\mu\\
 & =\alpha(\alpha-1)\int_{M}u^{\alpha-2}(t)(Lu(t))^{2}\xi_{k}d\mu+\alpha\int_{M}u^{\alpha-1}(t)LLu(t)\xi_{k}d\mu\\
 & =\alpha(\alpha-1)\int_{M}u^{\alpha-2}(t)(Lu(t))^{2}\xi_{k}d\mu+\alpha\int_{M}L\left(u^{\alpha-1}(t)\xi_{k}\right)Lu(t)d\mu\\
 & =\alpha(\alpha-1)\int_{M}u^{\alpha-2}(t)(Lu(t))^{2}\xi_{k}d\mu\\
 & \quad+\alpha\int_{M}\left(Lu^{\alpha-1}(t)\xi_{k}+L\xi_{k}u^{\alpha-1}(t)+2\langle\nabla u^{\alpha-1}(t),\nabla\xi_{k}\rangle\right)Lu(t)d\mu\\
 & =2\alpha(\alpha-1)\int_{M}u^{\alpha-2}(t)(Lu(t))^{2}\xi_{k}d\mu+\alpha(\alpha-1)(\alpha-2)\int_{M}u^{\alpha-3}(t)|\nabla u(t)|^{2}Lu(t)\xi_{k}d\mu\\
 & \quad+\alpha\int_{M}u^{\alpha-1}(t)Lu(t)L\xi_{k}d\mu+2\alpha(\alpha-1)\int_{M}u^{\alpha-2}(t)\langle\nabla u(t),\nabla\xi_{k}\rangle Lu(t)d\mu\\
 & =:2\alpha(\alpha-1)I_{3,k}(t)+\alpha(\alpha-1)(\alpha-2)I_{4,k}(t)+\alpha I_{5,k}(t)+2\alpha(\alpha-1)I_{6,k}(t).
\end{align*}
Thus, for any $0<a\le s<t\le b<\infty$,
\begin{align}
H'_{k}(t)-H'_{k}(s) & =2\alpha(\alpha-1)\int^{t}_{s}I_{3,k}(r)dr+\alpha(\alpha-1)(\alpha-2)\int^{t}_{s}I_{4,k}(r)dr\nonumber \\
 & \quad+\alpha\int^{t}_{s}I_{5,k}(r)dr+2\alpha(\alpha-1)\int^{t}_{s}I_{6,k}(r)dr.\label{eq:2-5-6}
\end{align}
By the monotone convergence theorem, it is easy to see that
\begin{equation}
I_{3,k}(t)\uparrow I_{3}(t):=\int_{M}u^{\alpha-2}(t)(Lu(t))^{2}d\mu.\label{eq:2-5-2}
\end{equation}
By Lemma \ref{thm:gradient estimates}, there exists constants $C_{1}(K,m,n,b),C_{2}(K,m,n,b)>0$,
such that
\begin{equation}
|Lu(x,t)|\le C_{1}(K,m,n,b)\left(\frac{d(x,y)}{t}+\frac{1}{\sqrt{t}}\right)^{2}u(x,t),\label{eq:3-12-13}
\end{equation}
\begin{equation}
|\nabla Lu(x,t)|\le C_{2}(K,m,n,b)\left(\frac{d(x,y)}{t}+\frac{1}{\sqrt{t}}\right)^{3}u(x,t).\label{eq:3-12-14}
\end{equation}
Substituting the estimate \eqref{eq:3-12-13}, we have
\begin{align*}
I_{3}(t) & =\int_{M}u^{\alpha-2}(t)(Lu(t))^{2}d\mu\\
 & \le C_{1}(K,m,n,b)^{2}\int_{M}u^{\alpha}(t)\left(\frac{d(x,y)}{t}+\frac{1}{\sqrt{t}}\right)^{4}d\mu\\
 & \le C_{1}(K,m,n,b)^{2}\sup_{a\le t\le b}\Vert u(\cdot,t)\Vert^{\alpha-1}_{\infty}\int_{M}u(t)\left(\frac{d(x,y)}{t}+\frac{1}{\sqrt{t}}\right)^{4}d\mu.
\end{align*}
Then it follows from \eqref{eq:3-12-1} and Lemma \ref{lem:Let-.-Let}
that
\[
\int_{M}u^{\alpha-2}(t)(Lu(t))^{2}d\mu<\infty,\quad\forall t>0,
\]
and
\begin{equation}
\sup_{a\le t\le b}I_{3}(t)<\infty.\label{eq:3-12-17}
\end{equation}
By \eqref{eq:2-5-2} and \eqref{eq:3-12-17}, we have
\begin{equation}
\int^{t}_{s}I_{3,k}(r)dr\to\int^{t}_{s}I_{3}(r)dr\label{eq:2-5-7}
\end{equation}
as $k\to\infty$, and
\begin{equation}
I_{3}(t)\in L^{1}([a,b]).\label{eq:2-7-1-1}
\end{equation}
Let $I_{4}(t):=\int_{M}u^{\alpha-3}(t)|\nabla u(t)|^{2}Lu(t)d\mu$.
 Substituting the estimates \eqref{eq:2-5-3} and \eqref{eq:3-12-13},
we have
\begin{align*}
|I_{4}(t)| & \le\int_{M}u^{\alpha-3}(t)|\nabla u(t)|^{2}|Lu(t)|d\mu\\
 & \le C^{2}_{0}(K,m,n,b)C_{1}(K,m,n,b)\int_{M}u^{\alpha-1}(t)\left(\frac{d(x,y)}{t}+\frac{1}{\sqrt{t}}\right)^{4}u(t)d\mu\\
 & \le C^{2}_{0}(K,m,n,b)C_{1}(K,m,n,b)\sup_{t\in[a,b]}\Vert u(\cdot,t)\Vert^{\alpha-1}_{\infty}\sup_{t\in[a,b]}\int_{M}\left(\frac{d(x,y)}{t}+\frac{1}{\sqrt{t}}\right)^{4}u(t)d\mu.
\end{align*}
Then it follows from \eqref{eq:3-12-1} and Lemma \ref{lem:Let-.-Let}
that
\[
\int_{M}u^{\alpha-3}(t)|\nabla u(t)|^{2}|Lu(t)|d\mu<\infty,\quad\forall t>0,
\]
and
\begin{equation}
\sup_{a\le t\le b}|I_{4}(t)|<\infty.\label{eq:2-7-1-4}
\end{equation}
Therefore,
\begin{equation}
I_{4}(t)\in L^{1}([a,b]),\label{eq:2-7-1-2}
\end{equation}
and by dominated convergence theorem,
\begin{equation}
\int^{b}_{a}I_{4,k}(r)dr\to\int^{b}_{a}I_{4}(r)dr\label{eq:2-5-8}
\end{equation}
as $k\to\infty$. 

Applying integration by parts,
\begin{align*}
I_{5,k} & =\int_{M}u^{\alpha-1}LuL\xi_{k}d\mu\\
 & =-\int_{M}\langle\nabla(u^{\alpha-1}Lu),\nabla\xi_{k}\rangle d\mu\\
 & =-(\alpha-1)\int_{M}u^{\alpha-2}Lu\langle\nabla u,\nabla\xi_{k}\rangle d\mu-\int_{M}u^{\alpha-1}\langle\nabla Lu,\nabla\xi_{k}\rangle d\mu\\
 & =:-(\alpha-1)J_{1,k}-J_{2,k}.
\end{align*}
By Hölder inequality, we have
\begin{align}
|J_{1,k}| & =\left|\int_{M}u^{\alpha-2}Lu\langle\nabla u,\nabla\xi_{k}\rangle d\mu\right|\nonumber \\
 & \le\frac{1}{k}\int_{M}u^{\alpha-2}|Lu||\nabla u|d\mu\nonumber \\
 & \le\frac{1}{k}\left(\int_{M}u^{\alpha-2}|Lu|^{2}d\mu\right)^{1/2}\left(\int_{M}u^{\alpha-2}|\nabla u|^{2}d\mu\right)^{1/2}.\label{eq:2-5-4}
\end{align}
From \eqref{eq:3-12-3} and \eqref{eq:3-12-17}, taking limits, we
derive
\begin{equation}
J_{1,k}\to0,\quad\int^{b}_{a}J_{1,k}(r)dr\to0,\label{eq:2-5-5}
\end{equation}
as $k\to\infty$.  

Substituting the estimate \eqref{eq:3-12-14},
\begin{align*}
|J_{2,k}| & =\left|\int_{M}u^{\alpha-1}\langle\nabla Lu,\nabla\xi_{k}\rangle d\mu\right|\\
 & \le\frac{1}{k}\int_{M}u^{\alpha-1}|\nabla Lu|d\mu\\
 & \le\frac{1}{k}C_{2}(K,m,n,b)\int_{M}u^{\alpha-1}\left(\frac{d(x,y)}{t}+\frac{1}{\sqrt{t}}\right)^{3}u(x,t)d\mu\\
 & \le\frac{1}{k}C_{2}(K,m,n,b)\sup_{t\in[a,b]}\Vert u(\cdot,t)\Vert^{\alpha-1}_{\infty}\int_{M}\left(\frac{d(x,y)}{t}+\frac{1}{\sqrt{t}}\right)^{3}u(x,t)d\mu.
\end{align*}
From \eqref{eq:3-12-1} and Lemma \ref{lem:Let-.-Let}, taking limits,
we derive
\[
J_{2,k}\to0,\quad\int^{b}_{a}J_{2,k}(r)dr\to0,
\]
 thus
\[
I_{5,k}\to0,
\]
and
\begin{equation}
\int^{b}_{a}I_{5,k}(r)dr\to0\label{eq:2-5-9}
\end{equation}
as $k\to\infty$.

Note that
\[
I_{6,k}=J_{1,k},
\]
then, by \eqref{eq:2-5-5},  we obtain
\begin{equation}
\int^{b}_{a}I_{6,k}(r)dr\to0\label{eq:2-5-10}
\end{equation}
as $k\to\infty$.

Letting $k\to\infty$ on both sides of \eqref{eq:2-5-6}, and then
using \eqref{eq:2-5-11}\eqref{eq:2-5-7}\eqref{eq:2-5-8}\eqref{eq:2-5-9}
and \eqref{eq:2-5-10}, we obtain
\begin{equation}
H'(t)-H'(s)=2\alpha(\alpha-1)\int^{t}_{s}I_{3}(r)dr+\alpha(\alpha-1)(\alpha-2)\int^{t}_{s}I_{4}(r)dr.\label{eq:2-7-1-3}
\end{equation}
Taking derivative with respect to $t$ on both sides of \eqref{eq:2-7-1-3},
and note that \eqref{eq:2-7-1-1}\eqref{eq:2-7-1-2} and $a<b$ are
arbitrary in $(0,\infty)$, we derive that
\[
H''(t)=2\alpha(\alpha-1)I_{3}(t)+\alpha(\alpha-1)(\alpha-2)I_{4}(t)
\]
holds for almost every $t>0$. Equivalently, \eqref{eq:2-9-2} holds
for almost every $t>0$. By an argument similar to that used to prove
the continuity of $I_{1}(t)$ on $(0,\infty)$, one readily verifies
that both $I_{3}(t)$ and $I_{4}(t)$ are continuous on $(0,\infty)$.
 It then follows by continuity that \eqref{eq:2-9-2} holds for every
$t>0$.
\end{IEEEproof}
\subsection{Proof for general positive solutions}

The preceding proof established the dissipation identities of Theorem~\ref{thm:derivatives of TsallisEnt-1} for fundamental solutions. For the positive solution of the heat equation, the corresponding
Tsallis entropy dissipation formulas are given as follows.
\begin{IEEEproof}[Proof of Theorem~\ref{thm:Let--be}]
The proof follows the same lines as that of Theorem \ref{thm:derivatives of TsallisEnt-1},
with the role of the fundamental solution replaced by that of a positive
solution. We first explain why the simplified assumption
\eqref{eq:2-9-2-4} is sufficient for all the cutoff and limiting
arguments used in that proof. Fix $0<a<b<\infty$ and write
\[
\begin{aligned}
U&:=\sup_{a\le t\le b}\int_Mu^\alpha\,d\mu,\qquad
B:=\int_a^b\!\int_Mu^{\alpha-2}|Lu|^2\,d\mu\,dt,\\
C&:=\int_a^b\!\int_Mu^{\alpha-2}|\nabla Lu|^2\,d\mu\,dt,\qquad
D:=\int_a^b\!\int_Mu^{\alpha-4}|\nabla u|^4\,d\mu\,dt.
\end{aligned}
\]
All four quantities are finite by \eqref{eq:2-9-2-3}
 and \eqref{eq:2-9-2-4}. By the Cauchy--Schwarz inequality,
\[
\begin{aligned}
\int_a^b\!\int_Mu^{\alpha-2}|\nabla u|^2\,d\mu\,dt
&\le\left(\int_a^b\!\int_Mu^\alpha\,d\mu\,dt\right)^{1/2}D^{1/2}\\
&\le\sqrt{(b-a)UD}<\infty.
\end{aligned}
\]
Similarly,
\[
\begin{aligned}
\int_a^b\!\int_Mu^{\alpha-3}|Lu||\nabla u|^2\,d\mu\,dt
&\le B^{1/2}D^{1/2}\\
&\le\frac12(B+D)<\infty.
\end{aligned}
\]
Consequently, \eqref{eq:2-9-2-4} together with
\eqref{eq:2-9-2-3} implies the original four-term
integrability bound used in the cutoff argument; in particular, the
mixed term needs no separate assumption. For the fundamental solution,
the same quartic-gradient bound follows directly from the heat-kernel
gradient estimate \eqref{eq:2-5-3}, the local $L^\infty$ estimate
\eqref{eq:3-12-1}, and the polynomial moment estimates in Lemma
\ref{lem:Let-.-Let}. These estimates ensure the necessary convergence
of the integrals and limits.

For completeness, the endpoint $\alpha=1$ can be treated by
writing $S_1(u)=-\int_Mu\log u\,d\mu$ and using a cutoff
function $\xi_k$. The additional logarithmic-moment assumption
implies $S_1(u)$ is finite and, by Cauchy--Schwarz, the first
boundary error is bounded by
\[
\frac1k\left(\int_a^b\!\int_Mu(\log u)^2\,d\mu\,dt\right)^{1/2}
\left(\int_a^b\!\int_M\frac{|\nabla u|^2}{u}\,d\mu\,dt\right)^{1/2}
\longrightarrow0.
\]
The second boundary error is controlled by the product of
$\int_a^b\!\int_M(Lu)^2/u$ and
$\int_a^b\!\int_M|\nabla u|^2/u$, which are finite by the
estimates above. The cutoff argument therefore yields, for
almost every $t>0$,
\[
\frac{d}{dt}S_1(u)=\int_M\frac{|\nabla u|^2}{u}\,d\mu,
\]
and
\[
\frac{d^2}{dt^2}S_1(u)
=-2\int_M\frac{(Lu)^2}{u}\,d\mu
+\int_M\frac{Lu|\nabla u|^2}{u^2}\,d\mu,
\]
which agree with \eqref{eq:3-12-15-1} and
\eqref{eq:2-9-2-2} at $\alpha=1$.
As in the proof of Theorem \ref{thm:derivatives of TsallisEnt-1},
the convergence of the integrals in \eqref{eq:3-12-15-1}
and \eqref{eq:2-9-2-2}, together with their continuity with respect to
$t$ on $(0,\infty)$, guarantees that these identities hold for every
$t>0$ under the continuity assumption in the theorem. To avoid repetition, we omit the details of the
proof.
\end{IEEEproof}

\section{\label{sec:Preliminaries-for-the}Auxiliary estimates for the concavity proofs}

\subsection{The nonlinear power substitution}

To establish the concavity of the entropy along the heat flow, we
adopt the idea introduced in \cite{sun2026concavity} and set $w:=u^{\alpha/2}$,
$\alpha=\frac{2}{1+\delta}$. A direct computation yields
\begin{equation}
\nabla w=\frac{\alpha}{2}u^{\frac{\alpha-2}{2}}\nabla u,\label{eq:3-5}
\end{equation}
\begin{equation}
|\nabla w|^{2}=\frac{\alpha^{2}}{4}u^{\alpha-2}|\nabla u|^{2},\label{eq:3-2}
\end{equation}
\begin{equation}
Lw=\frac{\alpha(\alpha-2)}{4}u^{\frac{\alpha}{2}-2}|\nabla u|^{2}
+\frac{\alpha}{2}u^{\frac{\alpha}{2}-1}Lu.
\label{eq:3-3}
\end{equation}
Then, under the assumptions of Theorem \ref{thm:derivatives of TsallisEnt-1},
the first order derivative expression \eqref{eq:3-12-15} can be rewritten
as
\begin{equation}
\frac{d}{dt}S_{\alpha}(u)=2(1+\delta)\int_{M}|\nabla w|^{2}d\mu,\label{eq:3-6}
\end{equation}
and the second order derivative expression \eqref{eq:2-9-2} can be
rewritten as
\begin{equation}
\frac{d^{2}}{dt^{2}}S_{\alpha}(u)=-4(1+\delta)\left(\int_{M}(Lw)^{2}d\mu+\delta\int_{M}Lw\frac{|\nabla w|^{2}}{w}d\mu\right).\label{eq:3-4}
\end{equation}
The identities \eqref{eq:3-6} and \eqref{eq:3-4} follow by
substituting \eqref{eq:3-2} and \eqref{eq:3-3} into the dissipation
formulas; this is the standard power-transformation calculation used
in the Euclidean setting in \cite{sun2026concavity}.

\subsection{Integration-by-parts estimates needed for the proofs}

We also need the following two lemmas.

\begin{lem}
\label{lem:For-.-.}Under the assumption in Theorem \ref{thm:derivatives of TsallisEnt-1}.
For $\delta\in(-1,1]$, $|\nabla w|,|\nabla^{2}w|,Lw,|\nabla Lw|\in L^{2}(M,\mu)$,
and $L|\nabla w|^{2},|\nabla|\nabla w|^{2}|\in L^{1}(M,\mu)$.
\end{lem}
\begin{IEEEproof}
For convenience, we recall the identities \eqref{eq:3-5} and \eqref{eq:3-3},
together with two further identities derived from them.
\[
\nabla w=\frac{\alpha}{2}u^{\frac{\alpha}{2}-1}\nabla u,
\]
\[
\nabla^{2}w=\frac{\alpha}{2}u^{\frac{\alpha}{2}-1}\nabla^{2}u+\frac{\alpha}{2}\left(\frac{\alpha}{2}-1\right)u^{\frac{\alpha}{2}-2}\nabla u\otimes\nabla u,
\]
\[
Lw=\frac{\alpha}{2}u^{\frac{\alpha}{2}-1}Lu+\frac{\alpha}{2}\left(\frac{\alpha}{2}-1\right)u^{\frac{\alpha}{2}-2}|\nabla u|^{2},
\]
\begin{align*}
\nabla Lw & =\frac{\alpha}{2}u^{\frac{\alpha}{2}-1}\nabla Lu+\frac{\alpha}{2}\left(\frac{\alpha}{2}-1\right)u^{\frac{\alpha}{2}-2}Lu\nabla u\\
 & \quad+\frac{\alpha}{2}\left(\frac{\alpha}{2}-1\right)u^{\frac{\alpha}{2}-2}\nabla|\nabla u|^{2}\\
 & \quad+\frac{\alpha}{2}\left(\frac{\alpha}{2}-1\right)\left(\frac{\alpha}{2}-2\right)u^{\frac{\alpha}{2}-3}|\nabla u|^{2}\nabla u.
\end{align*}
By Lemmas \ref{thm:gradient estimates-1} and \ref{thm:gradient estimates},
abbreviate $\frac{d(x,y)}{t}+\frac{1}{\sqrt{t}}$ as $A$, and fix
any $T>0$, for any $t\in(0,T]$, there exist constants $C_{1},C_{2},C_{3},C_{4}$
depending on $K,m,n,T,$ such that
\[
|\nabla w|\le C_{1}Au^{\frac{\alpha}{2}},
\]
\[
|\nabla^{2}w|\le C_{2}A^{2}u^{\frac{\alpha}{2}},
\]
\[
|Lw|\le C_{3}A^{2}u^{\frac{\alpha}{2}},
\]
\[
|\nabla Lw|\le C_{4}A^{3}u^{\frac{\alpha}{2}},
\]
where in the last inequality, we also use the fact that
\begin{equation}
|\nabla|\nabla u|^{2}|\le2|\nabla^{2}u||\nabla u|.\label{eq:3-10}
\end{equation}

Thus, for any $t\in(0,T]$,
\begin{equation}
|\nabla w|^{2}\le C^{2}_{1}A^{2}u^{\alpha},\label{eq:3-13}
\end{equation}
\begin{equation}
|\nabla^{2}w|^{2}\le C^{2}_{2}A^{4}u^{\alpha},\label{eq:3-14}
\end{equation}
\begin{equation}
|Lw|^{2}\le C^{2}_{3}A^{4}u^{\alpha},\label{eq:3-15}
\end{equation}
\begin{equation}
|\nabla Lw|^{2}\le C^{2}_{4}A^{6}u^{\alpha}.\label{eq:3-16}
\end{equation}
By Lemma \ref{lem:Let-.-Let}, for any fixed $y\in M$ and $t\in(0,T]$,
\begin{equation}
\int_{M}d^{k}(x,y)p_{t}(x,y)d\mu<+\infty,\quad\forall k\ge0,\label{eq:3-11}
\end{equation}
 Since $0<a<b<\infty$ are arbitrary in \eqref{eq:3-12-1}, it follows
that
\begin{equation}
\Vert u(\cdot,t)\Vert_{\infty}<+\infty\label{eq:3-12}
\end{equation}
for every $t>0$. Thereby, for any $k\ge0$,
\begin{align}
\int A^{k}u^{\alpha}d\mu & \le\int A^{k}u^{\alpha-1}ud\mu\nonumber \\
 & \le\Vert u\Vert^{\alpha-1}_{\infty}\int A^{k}ud\mu\nonumber \\
 & <+\infty,\label{eq:3-17}
\end{align}
where the last inequality comes from \eqref{eq:3-11} and \eqref{eq:3-12}.
Then $|\nabla w|,|\nabla^{2}w|,Lw$, and $|\nabla Lw|\in L^{2}(M,\mu)$ follow
from \eqref{eq:3-17} and estimates \eqref{eq:3-13}--\eqref{eq:3-16}.

By generalized Bochner--Weitzenböck formula,
\[
L|\nabla w|^{2}=2\langle\nabla w,\nabla Lw\rangle+2|\nabla^{2}w|^{2}+2Ric(L)(\nabla w,\nabla w).
\]
By the bounded geometry and $\nabla\phi\in C^{3}_{b}$,
\[
\Vert Ric(L)\Vert_{\infty}=\Vert Ric+\nabla^{2}\phi\Vert_{\infty}<\infty.
\]
Therefore, using the conclusions $|\nabla w|,|\nabla^{2}w|$, and
$|\nabla Lw|\in L^{2}(M,\mu)$ established above in this lemma, together
with the Cauchy--Schwarz inequality, we obtain
\begin{align*}
\int_{M}|L|\nabla w|^{2}|d\mu & \le2\Vert\nabla w\Vert_{2}\Vert\nabla Lw\Vert_{2}+2\Vert\nabla^{2}w\Vert^{2}_{2}\\
 & \quad+2\Vert Ric(L)\Vert_{\infty}\Vert\nabla w\Vert^{2}_{2}<\infty,
\end{align*}
that is, $L|\nabla w|^{2}\in L^{1}(M,\mu)$.

Finally, $|\nabla|\nabla w|^{2}|\in L^{1}(M,\mu)$ follows from $|\nabla|\nabla w|^{2}|\le2|\nabla^{2}w||\nabla w|$,
$|\nabla w|,|\nabla^{2}w|\in L^{2}(M,\mu)$, and Cauchy--Schwarz
inequality.
\end{IEEEproof}
\begin{lem}
\label{lem:For-.}Under the assumption in Theorem \ref{thm:derivatives of TsallisEnt-1}.
For $\delta\in(-1,1]$, the following integrability properties hold:
\[
|\nabla\sqrt{w}|^{2},\quad |\nabla|\nabla\sqrt{w}|^{2}|
\in L^{2}(M,\mu),
\]
and
\[
w|\nabla|\nabla\sqrt{w}|^{2}|,\quad
wL|\nabla\sqrt{w}|^{2}\in L^{1}(M,\mu).
\]
\end{lem}
\begin{IEEEproof}
It is straightforward to compute that
\begin{equation}
|\nabla\sqrt{w}|^{2}=\frac{1}{4}\frac{|\nabla w|^{2}}{w},\label{eq:3-18}
\end{equation}
\begin{equation}
\nabla|\nabla\sqrt{w}|^{2}=\frac{1}{4}\frac{\nabla|\nabla w|^{2}}{w}-\frac{1}{4}\frac{|\nabla w|^{2}\nabla w}{w^{2}},\label{eq:3-19}
\end{equation}
\begin{equation}
L|\nabla\sqrt{w}|^{2}=-\frac{1}{4}\frac{Lw|\nabla w|^{2}}{w^{2}}+\frac{1}{2}\frac{|\nabla w|^{4}}{w^{3}}+\frac{1}{4}\frac{L|\nabla w|^{2}}{w}-\frac{1}{2}\frac{\langle\nabla|\nabla w|^{2},\nabla w\rangle}{w^{2}}.\label{eq:3-20}
\end{equation}
 Using the same notations in the proof of Lemma \ref{lem:For-.-.}.
Note that $w=u^{\frac{\alpha}{2}}$. By \eqref{eq:3-13} and \eqref{eq:3-18},
\[
|\nabla w|^{4}\le C^{4}_{1}A^{4}u^{2\alpha},
\]
\[
|\nabla\sqrt{w}|^{4}=\frac{1}{16}\frac{|\nabla w|^{4}}{w^{2}}\le\frac{1}{16}C^{4}_{1}A^{4}u^{\alpha}.
\]
Thus, by \eqref{eq:3-17}, $|\nabla\sqrt{w}|^{2}\in L^{2}(M,\mu)$
follows. 

By \eqref{eq:3-19}\eqref{eq:3-13}\eqref{eq:3-14} and Cauchy-Schwarz
inequality,
\begin{align*}
\left|\nabla|\nabla\sqrt{w}|^{2}\right|^{2} & =\frac{1}{16}\frac{\left|\nabla|\nabla w|^{2}\right|^{2}}{w^{2}}+\frac{1}{16}\frac{|\nabla w|^{6}}{w^{4}}-\frac{1}{8}\frac{\langle\nabla|\nabla w|^{2},\nabla w\rangle|\nabla w|^{2}}{w^{3}}\\
 & \le\frac{1}{4}\frac{|\nabla^{2}w|^{2}|\nabla w|^{2}}{w^{2}}+\frac{1}{16}\frac{|\nabla w|^{6}}{w^{4}}+\frac{1}{4}\frac{|\nabla^{2}w||\nabla w|^{4}}{w^{3}}\\
 & \le\frac{1}{4}C^{2}_{1}C^{2}_{2}A^{6}u^{\alpha}+\frac{1}{16}C^{6}_{1}A^{6}u^{\alpha}+\frac{1}{4}C^{4}_{1}C_{2}A^{6}u^{\alpha},
\end{align*}
and
\begin{align*}
w|\nabla|\nabla\sqrt{w}|^{2}| & =\left|\frac{1}{4}\nabla|\nabla w|^{2}-\frac{1}{4}\frac{|\nabla w|^{2}\nabla w}{w}\right|\\
 & \le\frac{1}{2}|\nabla^{2}w||\nabla w|+\frac{1}{4}\frac{|\nabla w|^{3}}{w}\\
 & \le\frac{1}{2}C_{1}C_{2}A^{3}u^{\alpha}+\frac{1}{4}C^{3}_{1}A^{3}u^{\alpha}.
\end{align*}
Thus, by \eqref{eq:3-17}, $|\nabla|\nabla\sqrt{w}|^{2}|\in L^{2}(M,\mu)$
and $w|\nabla|\nabla\sqrt{w}|^{2}|\in L^{1}(M,\mu)$ follow. 

By \eqref{eq:3-20}, \eqref{eq:3-13}, \eqref{eq:3-14}, and \eqref{eq:3-15}
and Cauchy-Schwarz inequality,
\begin{align*}
w\left|L|\nabla\sqrt{w}|^{2}\right| & \le\frac{1}{4}\frac{|Lw||\nabla w|^{2}}{w}+\frac{1}{2}\frac{|\nabla w|^{4}}{w^{2}}+\frac{1}{4}\left|L|\nabla w|^{2}\right|+\frac{1}{2}\frac{\left|\nabla|\nabla w|^{2}\right||\nabla w|}{w}\\
 & \le\frac{1}{4}C^{2}_{1}C_{3}A^{4}u^{\alpha}+\frac{1}{2}C^{4}_{1}A^{4}u^{\alpha}+C^{2}_{1}C_{2}A^{4}u^{\alpha}+\frac{1}{4}\left|L|\nabla w|^{2}\right|.
\end{align*}
 by \eqref{eq:3-17} and Lemma \ref{lem:For-.-.}, $wL|\nabla\sqrt{w}|^{2}\in L^{1}(M,\mu)$
follows.
\end{IEEEproof}
Then we obtain the following useful integration by parts formula.
\begin{prop}
\label{prop:For-.-we}Under the assumption in Theorem \ref{thm:derivatives of TsallisEnt-1}.
For $\delta\in(-1,1]$, we have
\begin{align}
\int_{M}wLwd\mu & =-\int_{M}|\nabla w|^{2}d\mu,\label{eq:4-1-1}\\
\int_{M}L|\nabla w|^{2}d\mu & =0,\label{eq:4-1-2}\\
\int_{M}\langle\nabla Lw,\nabla w\rangle d\mu & =-\int_{M}|Lw|^{2}d\mu,\label{eq:4-1-3}\\
\int_{M}Lw|\nabla\sqrt{w}|^{2}d\mu & =-\int_{M}\langle\nabla w,\nabla|\nabla\sqrt{w}|^{2}\rangle d\mu\nonumber \\
 & =\int_{M}wL|\nabla\sqrt{w}|^{2}d\mu.\label{eq:4-1-4}
\end{align}
\end{prop}
\begin{IEEEproof}
As in the proof of Theorem \ref{thm:derivatives of TsallisEnt-1},
let $\{\xi_{k}\}^{\infty}_{k=1}$ be an increasing sequence of functions
in $C^{\infty}_{0}(M)$ such that $0\le\xi_{k}\le1$ and $|\nabla\xi_{k}|\le\frac{1}{k}$.
 
\begin{align*}
\int_{M}wLw\xi_{k}d\mu & =-\int_{M}\langle\nabla w,\nabla(w\xi_{k})\rangle d\mu\\
 & =-\int_{M}|\nabla w|^{2}\xi_{k}d\mu-\int_{M}\langle\nabla w,\nabla\xi_{k}\rangle wd\mu.
\end{align*}
By monotone convergence theorem,
\begin{equation}
\int_{M}|\nabla w|^{2}\xi_{k}d\mu\to\int_{M}|\nabla w|^{2}d\mu\label{eq:4-1-7}
\end{equation}
as $k\to\infty$. By dominated convergence theorem, 
\begin{equation}
\int_{M}wLw\xi_{k}d\mu\to\int_{M}wLwd\mu\label{eq:4-1-8}
\end{equation}
as $k\to\infty$ provided that
\[
\int_{M}w|Lw|d\mu<+\infty.
\]
Since
\begin{align*}
|\int_{M}\langle\nabla w,\nabla\xi_{k}\rangle wd\mu| & \le\frac{1}{k}\int_{M}|\nabla w|wd\mu,
\end{align*}
we find that
\begin{equation}
\int_{M}\langle\nabla w,\nabla\xi_{k}\rangle wd\mu\to0\label{eq:4-1-9}
\end{equation}
as $k\to\infty$ provided that
\[
\int_{M}|\nabla w|wd\mu<+\infty.
\]
 Indeed, by \eqref{eq:3-12-4} in the proof of Theorem \ref{thm:derivatives of TsallisEnt-1},
\begin{equation}
\Vert w\Vert_{2}=\left(\int_{M}u^{\alpha}d\mu\right)^{\frac{1}{2}}<+\infty.\label{eq:4-1-5}
\end{equation}
Then, by \eqref{eq:4-1-5} and Lemma \ref{lem:For-.-.}, and the
Cauchy-Schwarz inequality,
\[
\int_{M}w|Lw|d\mu\le\Vert w\Vert_{2}\Vert Lw\Vert_{2}<+\infty,
\]
and
\[
\int_{M}|\nabla w|wd\mu\le\Vert\nabla w\Vert_{2}\Vert w\Vert_{2}<+\infty.
\]
 Therefore, \eqref{eq:4-1-8} and \eqref{eq:4-1-9} hold and, together
with \eqref{eq:4-1-7}, imply \eqref{eq:4-1-1}.

Along the same lines, applying integrations by parts, we obtain
\[
\int_{M}L|\nabla w|^{2}\xi_{k}d\mu=-\int_{M}\langle\nabla|\nabla w|^{2},\nabla\xi_{k}\rangle d\mu.
\]
Thus, identity \eqref{eq:4-1-2} follows once we verify that
\[
\int_{M}|L|\nabla w|^{2}|d\mu<+\infty,
\]
and
\begin{equation}
\int_{M}|\nabla|\nabla w|^{2}|d\mu<+\infty.\label{eq:4-1-6}
\end{equation}
Lemma \ref{lem:For-.-.} yields the two inequalities above, thereby
proving \eqref{eq:4-1-2}.
\[
\]

Applying integrations by parts, we obtain
\begin{align*}
\int_{M}|Lw|^{2}\xi_{k}d\mu & =-\int_{M}\langle\nabla w,\nabla(\xi_{k}Lw)\rangle d\mu\\
 & =-\int_{M}\langle\nabla w,\nabla Lw\rangle\xi_{k}d\mu-\int_{M}\langle\nabla w,\nabla\xi_{k}\rangle Lwd\mu.
\end{align*}
Thus,  identity \eqref{eq:4-1-3} follows once we verify that
\[
\int_{M}|\nabla w||\nabla Lw|d\mu<+\infty,
\]
and
\[
\int_{M}|\nabla w||Lw|d\mu<+\infty.
\]
By Lemma \ref{lem:For-.-.}, and the Cauchy-Schwarz inequality,
\[
\int_{M}|\nabla w||\nabla Lw|d\mu\le\Vert\nabla w\Vert_{2}\Vert\nabla Lw\Vert_{2}<+\infty,
\]
\[
\int_{M}|\nabla w||Lw|d\mu\le\Vert\nabla w\Vert_{2}\Vert Lw\Vert_{2}<+\infty.
\]
 This proves \eqref{eq:4-1-3}.

Applying integrations by parts, we obtain
\begin{align*}
\int_{M}Lw|\nabla\sqrt{w}|^{2}\xi_{k}d\mu & =-\int_{M}\langle\nabla w,\nabla\left(|\nabla\sqrt{w}|^{2}\xi_{k}\right)\rangle d\mu\\
 & =-\int_{M}\langle\nabla w,\nabla|\nabla\sqrt{w}|^{2}\rangle\xi_{k}d\mu-\int_{M}\langle\nabla w,\nabla\xi_{k}\rangle|\nabla\sqrt{w}|^{2}d\mu,
\end{align*}
\begin{align*}
\int_{M}wL|\nabla\sqrt{w}|^{2}\xi_{k}d\mu & =-\int_{M}\langle\nabla(w\xi_{k}),\nabla|\nabla\sqrt{w}|^{2}\rangle d\mu\\
 & =-\int_{M}\langle\nabla w,\nabla|\nabla\sqrt{w}|^{2}\rangle\xi_{k}d\mu-\int_{M}\langle\nabla\xi_{k},\nabla|\nabla\sqrt{w}|^{2}\rangle wd\mu.
\end{align*}
Thus,  identity \eqref{eq:4-1-4} follows once we verify that
\[
\int_{M}|Lw||\nabla\sqrt{w}|^{2}d\mu<+\infty,
\]
\[
\int_{M}|\nabla w||\nabla|\nabla\sqrt{w}|^{2}|d\mu<+\infty,
\]
\[
\int_{M}|\nabla w||\nabla\sqrt{w}|^{2}d\mu<+\infty,
\]
\[
\int_{M}w|L|\nabla\sqrt{w}|^{2}|d\mu<+\infty,
\]
\[
\int_{M}|\nabla|\nabla\sqrt{w}|^{2}|wd\mu<+\infty.
\]
Lemmas \ref{lem:For-.-.} and \ref{lem:For-.}, and Cauchy-Schwarz
inequality yield these five inequalities above, thereby proving \eqref{eq:4-1-4}.
\end{IEEEproof}
\subsection{Curvature-dependent inequalities for the main proofs}

The inequalities in the following two propositions generalize those
in Propositions 14 and 16 in \cite{sun2026concavity} to the manifold
setting, and their proofs are quite similar.
\begin{prop}
\label{prop:Under-the-assumption-1}Under the assumption in Theorem
\ref{thm:derivatives of TsallisEnt-1}. For $\delta\in(-1,1]$,
\begin{equation}
\int_{M}|\nabla^{2}w|^{2}d\mu=\int_{M}|Lw|^{2}d\mu-\int_{M}\langle Ric(L)\nabla w,\nabla w\rangle d\mu,\label{eq:2-9-4}
\end{equation}
\begin{equation}
\int_{M}Lw\frac{|\nabla w|^{2}}{w}d\mu=\int_{M}\frac{|\nabla w|^{4}}{w^{2}}d\mu-2\int_{M}\frac{\nabla^{2}w(\nabla w,\nabla w)}{w}d\mu,\label{eq:2-9-5}
\end{equation}
\begin{equation}
\int_{M}|Lw|^{2}d\mu+\int_{M}Lw\frac{|\nabla w|^{2}}{w}d\mu=\int_{M}w^{2}|\nabla^{2}\log w|^{2}d\mu+\int_{M}\langle Ric(L)\nabla w,\nabla w\rangle d\mu.\label{eq:2-9-5-1}
\end{equation}
Furthermore, suppose that
\[
Ric(L)\ge0,
\]
then
\begin{align}
\int_{M}|Lw|^{2}d\mu+\int_{M}Lw\frac{|\nabla w|^{2}}{w}d\mu & \ge\int_{M}|\nabla^{2}w-\frac{\nabla w\otimes\nabla w}{w}|^{2}d\mu\nonumber \\
 & =\int_{M}w^{2}|\nabla^{2}\log w|^{2}d\mu\ge0.\label{eq:5-2}
\end{align}
\end{prop}
\begin{IEEEproof}
By the generalized Bochner--Weitzenböck formula
\[
\frac{1}{2}L|\nabla w|^{2}=\langle\nabla w,\nabla Lw\rangle+|\nabla^{2}w|^{2}+\langle Ric(L)\nabla w,\nabla w\rangle,
\]
together with the integration by parts formulas \eqref{eq:4-1-2} and \eqref{eq:4-1-3}
in Proposition \ref{prop:For-.-we}, we have
\begin{align*}
\int_{M}|\nabla^{2}w|^{2}d\mu & =\frac{1}{2}\int_{M}L|\nabla w|^{2}d\mu-\int_{M}\langle\nabla w,\nabla Lw\rangle d\mu-\int_{M}\langle Ric(L)\nabla w,\nabla w\rangle d\mu\\
 & =\int_{M}|Lw|^{2}d\mu-\int_{M}\langle Ric(L)\nabla w,\nabla w\rangle d\mu.
\end{align*}
By integrations by parts formula \eqref{eq:4-1-4} in Proposition
\ref{prop:For-.-we}, we have
\begin{align*}
\int_{M}Lw\frac{|\nabla w|^{2}}{w}d\mu & =-\int_{M}\left\langle \nabla w,\nabla\left(\frac{|\nabla w|^{2}}{w}\right)\right\rangle d\mu\\
 & =-2\int_{M}\frac{\nabla^{2}w(\nabla w,\nabla w)}{w}d\mu+\int_{M}\frac{|\nabla w|^{4}}{w^{2}}d\mu.
\end{align*}
Thus,
\begin{align*}
\int_{M}|Lw|^{2}d\mu+\int_{M}Lw\frac{|\nabla w|^{2}}{w}d\mu & =\int_{M}|\nabla^{2}w|^{2}d\mu+\int_{M}\langle Ric(L)\nabla w,\nabla w\rangle d\mu\\
 & \quad-2\int_{M}\frac{\nabla^{2}w(\nabla w,\nabla w)}{w}d\mu+\int_{M}\frac{|\nabla w|^{4}}{w^{2}}d\mu\\
 & =\int_{M}|\nabla^{2}w-\frac{\nabla w\otimes\nabla w}{w}|^{2}d\mu+\int_{M}\langle Ric(L)\nabla w,\nabla w\rangle d\mu.
\end{align*}
Therefore,
\[
\int_{M}|Lw|^{2}d\mu+\int_{M}Lw\frac{|\nabla w|^{2}}{w}d\mu\ge\int_{M}|\nabla^{2}w-\frac{\nabla w\otimes\nabla w}{w}|^{2}d\mu\ge0
\]
if
\[
Ric(L)\ge0.
\]
Since 
\[
\nabla^{2}w-\frac{\nabla w\otimes\nabla w}{w}=w\nabla^{2}\log w,
\]
the conclusion follows.
\end{IEEEproof}
\begin{prop}
\label{prop:Under-the-assumption}Under the assumption in Theorem
\ref{thm:derivatives of TsallisEnt-1}. For $\delta\in(-1,1]$,
\[
3\int_{M}(Lw)^{2}d\mu-\int_{M}Lw\frac{|\nabla w|^{2}}{w}d\mu\ge\int\langle Ric(L)\nabla w,\nabla w\rangle d\mu.
\]
Furthermore, suppose that
\[
Ric(L)\ge0,
\]
we have
\begin{equation}
3\int_{M}(Lw)^{2}d\mu\ge\int_{M}Lw\frac{|\nabla w|^{2}}{w}d\mu.\label{eq:5-3}
\end{equation}
\end{prop}
\begin{IEEEproof}
Let $q=m-n$, $c_{u}=\langle\nabla\phi,\nabla w\rangle$.
We first write the proof for a positive integer $q$. For arbitrary
real $q>0$, the same calculation is valid algebraically with the
auxiliary scalar quadratic form $c_u^2/q$, and no integer-dimensional
auxiliary space is required. When $q=0$, the weight is constant, and
the auxiliary block below is omitted. Define the
symmetric matrices in $TM\oplus\mathbb{R}^{q}$ 
\[
P:=\begin{pmatrix}\nabla^{2}w & 0\\
0 & -\frac{c_{u}}{q}I_{q}
\end{pmatrix},\quad Q:=\begin{pmatrix}\frac{\nabla w\otimes\nabla w}{w} & 0\\
0 & 0
\end{pmatrix},
\]
\[
P-\frac{Q}{3}=\begin{pmatrix}\nabla^{2}w-\frac{\nabla w\otimes\nabla w}{3w} & 0\\
0 & -\frac{c_{u}}{q}I_{q}
\end{pmatrix},
\]
\begin{equation}
\left|P-\frac{Q}{3}\right|^{2}=\left|\nabla^{2}w-\frac{\nabla w\otimes\nabla w}{3w}\right|^{2}+\frac{\langle\nabla\phi,\nabla w\rangle^{2}}{m-n},\label{eq:2-9-6}
\end{equation}
and their traces
\begin{align*}
a & :=\text{tr}P=\text{tr}(\nabla^{2}w)-q\frac{c_{u}}{q}=\Delta w-\langle\nabla\phi,\nabla w\rangle=Lw,\\
b & :=\text{tr}Q=\text{tr}\left(\frac{\nabla w\otimes\nabla w}{w}\right)=\frac{|\nabla w|^{2}}{w}.
\end{align*}
Moreover, we define
\[
P:Q=\langle P,Q\rangle=\langle\nabla^{2}w,\frac{\nabla w\otimes\nabla w}{w}\rangle=\frac{\nabla^{2}w(\nabla w,\nabla w)}{w}.
\]
Then the following identity holds,
\begin{align}
2\left(a-\frac{b}{3}\right)^{2}+\left|P-\frac{Q}{3}\right|^{2} & =3a^{2}-ab+(|P|^{2}-a^{2})+\frac{1}{3}(b^{2}-ab-2P:Q).\label{eq:2-9-3}
\end{align}
This is because, expanding the left hand side gives
\[
2a^{2}-\frac{4}{3}ab+\frac{2}{9}b^{2}+|P|^{2}-\frac{2}{3}P:Q+\frac{1}{9}|Q|^{2},
\]
and using $|Q|^{2}=b^{2}$ gives the right hand side of \eqref{eq:2-9-3}.

By the definition of $P$,
\begin{align*}
|P|^{2} & =|\nabla^{2}w|^{2}+q\left(\frac{c_{u}}{q}\right)^{2}\\
 & =|\nabla^{2}w|^{2}+\frac{\langle\nabla\phi,\nabla w\rangle^{2}}{m-n}.
\end{align*}
Thus, \eqref{eq:2-9-4} gives
\begin{align}
\int_{M}|P|^{2}-a^{2}d\mu & =\int_{M}|\nabla^{2}w|^{2}+\frac{\langle\nabla\phi,\nabla w\rangle^{2}}{m-n}d\mu-\int_{M}(Lw)^{2}d\mu\nonumber \\
 & =\int_{M}\frac{\langle\nabla\phi,\nabla w\rangle^{2}}{m-n}d\mu-\int_{M}\langle Ric(L)\nabla w,\nabla w\rangle d\mu\nonumber \\
 & =-\int_{M}\langle Ric_{m,n}(L)\nabla w,\nabla w\rangle d\mu.\label{eq:2-9-7}
\end{align}

By the definition of $a,b,P:Q$ and \eqref{eq:2-9-5},
\begin{equation}
\int_{M}(b^{2}-ab-2P:Q)d\mu=\int_{M}\frac{|\nabla w|^{4}}{w^{2}}d\mu-\int_{M}Lw\frac{|\nabla w|^{2}}{w}d\mu-2\int_{M}\frac{\nabla^{2}w(\nabla w,\nabla w)}{w}d\mu=0.\label{eq:2-9-8}
\end{equation}

Integrate both sides of \eqref{eq:2-9-3}, and by \eqref{eq:2-9-6}, \eqref{eq:2-9-7}, and \eqref{eq:2-9-8},
we have
\begin{align*}
 & \int_{M}2\left(Lw-\frac{|\nabla w|^{2}}{3w}\right)^{2}+\left|\nabla^{2}w-\frac{\nabla w\otimes\nabla w}{3w}\right|^{2}+\frac{\langle\nabla\phi,\nabla w\rangle^{2}}{m-n}d\mu\\
= & 3\int_{M}(Lw)^{2}d\mu-\int_{M}Lw\frac{|\nabla w|^{2}}{w}d\mu-\int\langle Ric_{m,n}(L)\nabla w,\nabla w\rangle d\mu.
\end{align*}
Removing $\int_{M}\frac{\langle\nabla\phi,\nabla w\rangle^{2}}{m-n}d\mu$
on both sides gives 
\begin{align*}
 & \int_{M}2\left(Lw-\frac{|\nabla w|^{2}}{3w}\right)^{2}+\left|\nabla^{2}w-\frac{\nabla w\otimes\nabla w}{3w}\right|^{2}d\mu\\
= & 3\int_{M}(Lw)^{2}d\mu-\int_{M}Lw\frac{|\nabla w|^{2}}{w}d\mu-\int\langle Ric(L)\nabla w,\nabla w\rangle d\mu.
\end{align*}
The left hand side is a sum of squares, hence nonnegative. Then we
complete the proof.
\end{IEEEproof}

\section{\label{sec:The-concavity-of}Proofs of the Tsallis entropy concavity theorems}

Under the assumptions of Theorem \ref{thm:derivatives of TsallisEnt-1},
it is clear that the Tsallis entropy is nondecreasing along the heat
flow. However, the sign of its second derivative in \eqref{eq:2-9-2}
cannot be determined directly and requires further investigation.
Based on the preliminaries in Section \ref{sec:Preliminaries-for-the},
we verify the concavity of Tsallis entropy for the heat equation on
complete Riemannian manifolds with suitable curvature-dimension condition.

\subsection{Proof for the fundamental solution}

\begin{IEEEproof}[Proof of Theorem~\ref{thm:Under-the-assumptions}]
Recall that, let $w=u^{\alpha/2},\alpha=\frac{2}{1+\delta}$, where
$\delta\in(-1,1]$, then $\frac{\partial^{2}}{\partial t^{2}}h_{\alpha}(u)$
has an equivalent expression \eqref{eq:3-4},
\[
\frac{d^{2}}{dt^{2}}S_{\alpha}(u)=-4(1+\delta)\left(\int_{M}(Lw)^{2}d\mu+\delta\int_{M}Lw\frac{|\nabla w|^{2}}{w}d\mu\right).
\]
Therefore, it is enough to prove that
\[
\int_{M}(Lw)^{2}d\mu+\delta\int_{M}Lw\frac{|\nabla w|^{2}}{w}d\mu\ge0.
\]
If $0\le\delta\le1$, and $Ric(L)\ge0,$ \eqref{eq:5-2} gives
\begin{align*}
 & \int_{M}(Lw)^{2}d\mu+\delta\int_{M}Lw\frac{|\nabla w|^{2}}{w}d\mu\\
= & (1-\delta)\int_{M}(Lw)^{2}d\mu+\delta\left(\int_{M}(Lw)^{2}d\mu+\int_{M}Lw\frac{|\nabla w|^{2}}{w}d\mu\right)\\
 & \ge0.
\end{align*}
If $-\frac{1}{3}\le\delta\le0$, and $Ric(L)\ge0,$ \eqref{eq:5-3}
gives
\begin{align*}
 & \int_{M}(Lw)^{2}d\mu+\delta\int_{M}Lw\frac{|\nabla w|^{2}}{w}d\mu\\
\ge & \int_{M}(Lw)^{2}d\mu+3\delta\int_{M}(Lw)^{2}d\mu\\
= & (1+3\delta)\int_{M}(Lw)^{2}d\mu\\
\ge & 0.
\end{align*}
Note that $\alpha=\frac{2}{1+\delta}$, thus, $0\le\delta\le1$ corresponds
to $1\le\alpha\le2$, $-\frac{1}{3}\le\delta\le0$ corresponds to
$2\le\alpha\le3.$ Then the conclusion follows.
\end{IEEEproof}
\subsection{Proof for general positive solutions}

Theorem~\ref{thm:Under-the-assumptions} covers the fundamental solution of the heat equation,
rather than for an arbitrary positive solution. In the Shannon entropy
setting, \cite{li2016hamilton} likewise establishes the entropy dissipation
formula and the entropic concavity result for the fundamental solution
of the heat equation. For the positive solution $u$ of the heat equation,
we also set $w=u^{\alpha/2}$, $\alpha=\frac{2}{1+\delta}$, the corresponding
entropic concavity result is given as follows.
\begin{IEEEproof}[Proof of Theorem~\ref{thm:Under-the-assumptions-2}]
For the fundamental solution $u$ of the heat equation, the proof
of Theorem \ref{thm:Under-the-assumptions} relies on rewriting the
second time derivative of Tsallis entropy in the form \eqref{eq:3-4}
via the change of variable $w=u^{\alpha/2}$, $\alpha=\frac{2}{1+\delta}$,
and then determining the sign of \eqref{eq:3-4} using the functional
inequalities \eqref{eq:5-2} in Proposition \ref{prop:Under-the-assumption-1}
and \eqref{eq:5-3} in Proposition \ref{prop:Under-the-assumption}.
Note also that \eqref{eq:5-2} and \eqref{eq:5-3} follow from the
integration by parts formulas \eqref{eq:4-1-2}, \eqref{eq:4-1-3}, and \eqref{eq:4-1-4}
in Proposition \ref{prop:For-.-we}. For a positive solution of the
heat equation, the same change of variable allows us to rewrite the
second time derivative of the Tsallis entropy in the same form. Therefore,
the proof of Theorem \ref{thm:Under-the-assumptions-2} follows along
the same lines as that of Theorem \ref{thm:Under-the-assumptions},
with the fundamental solution replaced by a positive solution of the
heat equation.
\end{IEEEproof}

\section{\label{sec:Tsallis-entropy-of}Proofs of the order-two complete-monotonicity theorems}

\subsection{Basic dissipation identities and an independent cutoff proof}

In this section, we restrict to the Tsallis entropy of order 2, that
is
\[
S_{2}(u)=1-\int_{M}u^{2}(x)d\mu(x).
\]

For the fundamental solution of the heat equation, the first two
time derivatives of $S_{2}$ are given as follows.
\begin{cor}
\label{cor:Let--be}Under the assumptions of Theorem \ref{thm:derivatives of TsallisEnt-1},
we have
\[
\frac{d}{dt}S_{2}(u)=2\int_{M}|\nabla u|^{2}d\mu,
\]
and
\begin{equation}
\frac{d^{2}}{dt^{2}}S_{2}(u)=-4\int_{M}(Lu)^{2}d\mu,\label{eq:2-9-2-1}
\end{equation}
for every $t>0$.
\end{cor}
\begin{IEEEproof}
This is a straightforward corollary of Theorem \ref{thm:derivatives of TsallisEnt-1}.
\end{IEEEproof}
For a positive solution of the heat equation,
the first two time derivatives of $S_{2}$ are given as follows.
\begin{cor}
\label{cor:order-two-positive-solution}
Let $(M,g)$ be an $n$-dimensional complete Riemannian manifold,
and $\phi\in C^{2}(M)$. Let $u$ be a smooth positive solution of
the heat equation
\[
\partial_{t}u=Lu.
\]
Assume that, for every $0<a<b<\infty$,
\begin{equation}
\sup_{t\in[a,b]}\int_{M}u^{2}d\mu<\infty,\label{eq:2-9-2-3-1}
\end{equation}
and
\begin{equation}
\int^{b}_{a}\int_{M}\left(|\nabla u|^{2}+|Lu|^{2}+|\nabla Lu|^{2}\right)d\mu dt<\infty.\label{eq:2-9-2-4-1}
\end{equation}
Then
\begin{equation}
\frac{d}{dt}S_{2}(u)=2\int_{M}|\nabla u|^{2}d\mu,\label{eq:2-9-2-5}
\end{equation}
and
\begin{equation}
\frac{d^{2}}{dt^{2}}S_{2}(u)=-4\int_{M}(Lu)^{2}d\mu,\label{eq:2-9-2-6}
\end{equation}
for almost every $t>0$. Furthermore, if the integrals in \eqref{eq:2-9-2-5} and \eqref{eq:2-9-2-6}
are convergent and continuous for every $t>0$, then \eqref{eq:2-9-2-5} and \eqref{eq:2-9-2-6}
hold for every $t>0$.
\end{cor}
\begin{IEEEproof}
We prove the identities directly, since the $L^2$ assumptions in this
corollary are weaker than the quartic-gradient assumption
\eqref{eq:2-9-2-4} for general $\alpha$. Fix $0<a<b<\infty$ and
choose the cutoff functions $\xi_R\in C_0^\infty(M)$ used in the
proof of Theorem \ref{thm:derivatives of TsallisEnt-1}, with
$0\le\xi_R\le1$, $\xi_R\uparrow1$, and
$|\nabla\xi_R|\le1/R$. Set
\[
H_R(t):=\int_Mu(t)^2\xi_R\,d\mu,\qquad
J_R(t):=\int_M|\nabla u(t)|^2\xi_R\,d\mu.
\]
Since $u_t=Lu$, integration by parts with compactly supported
$\xi_R$ gives
\[
\begin{aligned}
H_R'(t)&=2\int_MuLu\xi_R\,d\mu\\
&=-2J_R(t)-2\int_Mu\langle\nabla u,\nabla\xi_R\rangle\,d\mu.
\end{aligned}
\]
The last term satisfies
\[
\begin{aligned}
&\int_a^b\left|\int_Mu\langle\nabla u,\nabla\xi_R\rangle\,d\mu\right|dt\\
&\quad\le\frac1R
\left(\int_a^b\!\int_Mu^2\,d\mu\,dt\right)^{1/2}
\left(\int_a^b\!\int_M|\nabla u|^2\,d\mu\,dt\right)^{1/2}
\longrightarrow0.
\end{aligned}
\]
Passing to the limit in the integrated identity yields
$\frac d{dt}\int_Mu(t)^2d\mu=-2\int_M|\nabla u(t)|^2d\mu$
for almost every $t>0$, and hence \eqref{eq:2-9-2-5}.

For the second derivative, use $\partial_t\nabla u=\nabla Lu$
and integrate by parts once more:
\[
\begin{aligned}
J_R'(t)&=2\int_M\langle\nabla u,\nabla Lu\rangle\xi_R\,d\mu\\
&=-2\int_M(Lu)^2\xi_R\,d\mu
-2\int_MLu\langle\nabla u,\nabla\xi_R\rangle\,d\mu.
\end{aligned}
\]
By Cauchy--Schwarz and the second integrability assumption,
\[
\begin{aligned}
&\int_a^b\left|\int_MLu\langle\nabla u,\nabla\xi_R\rangle\,d\mu\right|dt\\
&\quad\le\frac1R
\left(\int_a^b\!\int_M|Lu|^2\,d\mu\,dt\right)^{1/2}
\left(\int_a^b\!\int_M|\nabla u|^2\,d\mu\,dt\right)^{1/2}
\longrightarrow0.
\end{aligned}
\]
Thus $\frac d{dt}\int_M|\nabla u(t)|^2d\mu
=-2\int_M(Lu(t))^2d\mu$ for almost every $t>0$.
Combining the two identities proves \eqref{eq:2-9-2-6}.
The stated continuity assumptions extend the almost-everywhere
identities to every positive time.
\end{IEEEproof}
\subsection{Higher derivatives via the self-adjoint generator}

For the fundamental solution of the heat equation, the higher order
time derivatives of $S_{2}$ are given as follows.
\begin{thm}
\label{thm:all derivative}Let $(M,g)$ be an $n$-dimensional complete
Riemannian manifold.  For any positive integer $k$, assume that
$\Vert\nabla^{j}Riem\Vert_{\infty}<C^{(1)}_{j}$ for some constant
$C^{(1)}_{j}>0,0\le j\le4k-1$. Let $\phi\in C^{4k}(M)$ be such that
$\nabla\phi\in C^{4k-1}_{b}(M)$, and suppose that there exist two
constants $m\ge n$ and $K\ge0$ such that
\[
Ric_{m,n}(L)\ge-K.
\]
Then, for any fundamental solution $u(x,t)$ to the heat equation
\[
\partial_{t}u=Lu,
\]
we have
\begin{equation}
\frac{d^{2k-1}}{dt^{2k-1}}S_{2}(u(t))=2^{2k-1}\int_{M}|\nabla L^{k-1}u|^{2}d\mu,\label{eq:7}
\end{equation}
\begin{equation}
\frac{d^{2k}}{dt^{2k}}S_{2}(u(t))=-2^{2k}\int_{M}(L^{k}u)^{2}d\mu,\label{eq:8}
\end{equation}
where $L^{0}$ denotes identity operator $\mathrm{id}$.
\end{thm}
\begin{IEEEproof}
Fix an arbitrary point $y\in M$, and let $u(t)=u(x,t)=p_{t}(x,y)$
be the fundamental solution. Let $0<a<b<\infty$ be arbitrary, and
fix $a\le s<t\le b$. 

For simplicity of presentation, we carry out the differentiations
under the integral sign and the integrations by parts directly. These
operations can be justified rigorously by the same cutoff argument
as that used in the proof of Theorem \ref{thm:derivatives of TsallisEnt-1},
together with repeated applications of Lemmas \ref{thm:gradient estimates-1},
\ref{thm:gradient estimates}, and \ref{lem:Let-.-Let}. More precisely,
one may introduce an increasing sequence $\{\xi_{m}\}^{\infty}_{m=1}\subset C^{\infty}_{0}(M)$
satisfying
\[
0\le\xi_{m}\le1,\quad\xi_{m}\uparrow1,\quad|\nabla\xi_{m}|\le\frac{1}{m},
\]
 and perform all computations after multiplying the relevant integrands
by $\xi_{m}$. Whenever a term involving $L\xi_{m}$ arises, the symmetry
of $L$ with respect to $d\mu$ allow $L$ to be transferred to the
other factor by integration by parts. On every fixed time interval
$[a,b]\subset(0,\infty)$, Lemmas \ref{thm:gradient estimates-1}
and \ref{thm:gradient estimates} provide the required high-order
derivative estimates for the heat kernel, while Lemma \ref{lem:Let-.-Let}
guarantees the uniform finiteness of its polynomial distance moments.
As in the proof of Theorem \ref{thm:derivatives of TsallisEnt-1},
these three lemmas imply the $L^{1}$-integrability of all terms arising
from the differentiations and integrations by parts. Consequently,
the error terms involving $\nabla\xi_{m}$ vanish as $m\to\infty$,
whereas the remaining terms converge by dominated convergence. We
may therefore omit the cutoff functions in what follows and carry
out the differentiations under the integral sign and the integrations
by parts directly.

By Corollary \ref{cor:Let--be}, suppose \eqref{eq:7} and \eqref{eq:8}
hold for $k$, then
\begin{align}
\frac{d^{2k+1}}{dt^{2k+1}}S_{2}(u(t)) & =\frac{d}{dt}\left(\frac{d^{2k}}{dt^{2k}}S_{2}(u(t))\right)\nonumber \\
 & =-2^{2k}\frac{d}{dt}\left(\int_{M}(L^{k}u)^{2}d\mu\right)\nonumber \\
 & =-2^{2k+1}\int_{M}L^{k}uL^{k+1}ud\mu,\label{eq:9}
\end{align}
and
\begin{align}
\frac{d^{2k+2}}{dt^{2k+2}}S_{2}(u(t)) & =\frac{d}{dt}\left(\frac{d^{2k+1}}{dt^{2k+1}}S_{2}(u(t))\right)\nonumber \\
 & =2^{2k+1}\frac{d}{dt}\left(\int_{M}|\nabla L^{k}u|^{2}d\mu\right)\nonumber \\
 & =2^{2k+2}\int_{M}\langle\nabla L^{k+1}u,\nabla L^{k}u\rangle d\mu.\label{eq:13}
\end{align}
By integration by parts,
\begin{align*}
\int_{M}L^{k}uL^{k+1}ud\mu & =\int_{M}L^{k}uL(L^{k}u)d\mu\\
 & =-\int_{M}|\nabla L^{k}u|^{2}d\mu,
\end{align*}
and
\[
\int_{M}(L^{k+1}u)^{2}d\mu=-\int_{M}\langle\nabla L^{k+1}u,\nabla L^{k}u\rangle d\mu.
\]
Thus,
\[
\frac{d^{2k+1}}{dt^{2k+1}}S_{2}(u(t))=2^{2k+1}\int_{M}|\nabla L^{k}u|^{2}d\mu,
\]
and
\[
\frac{d^{2k+2}}{dt^{2k+2}}S_{2}(u(t))=-2^{2k+2}\int_{M}(L^{k+1}u)^{2}d\mu.
\]
\textcolor{blue}{ }This completes the proof by induction.
\end{IEEEproof}
For the positive solution of the heat equation, the higher order time
derivatives of $S_{2}$ are given as follows.
\begin{thm}
\label{thm:Let--be-1}Let $(M,g)$ be an $n$-dimensional complete
Riemannian manifold, and $\phi\in C^{2}(M)$. Let $u$ be a smooth
positive solution of the heat equation
\[
\partial_{t}u=Lu.
\]
Let $k$ be an arbitrary positive integer. Assume that, for every $0<a<b<\infty$,
\begin{equation}
\sup_{t\in[a,b]}\int_{M}u^{2}d\mu<\infty,\label{eq:2-9-2-3-1-1}
\end{equation}
and
\begin{equation}
\int^{b}_{a}\int_{M}\left(|L^{i}u|^{2}+|\nabla L^{i}u|^{2}\right)d\mu dt<\infty\label{eq:2-9-2-4-1-1}
\end{equation}
holds for each $0\le i\le k$. Then
\begin{equation}
\frac{d^{2i-1}}{dt^{2i-1}}S_{2}(u(t))=2^{2i-1}\int_{M}|\nabla L^{i-1}u|^{2}d\mu\label{eq:2-9-2-5-1}
\end{equation}
and
\begin{equation}
\frac{d^{2i}}{dt^{2i}}S_{2}(u(t))=-2^{2i}\int_{M}(L^{i}u)^{2}d\mu,\label{eq:2-9-2-6-1}
\end{equation}
for each $1\le i\le k$ and almost every $t>0$. Furthermore, if the
integrals in \eqref{eq:2-9-2-5-1} and \eqref{eq:2-9-2-6-1} converge
and are continuous with respect to $t$ for each $1\le i\le k$ and
every $t>0$, then \eqref{eq:2-9-2-5-1} and \eqref{eq:2-9-2-6-1}
hold for each $1\le i\le k$ and every $t>0$.
\end{thm}
\begin{IEEEproof}
The proof follows the same lines as that of Theorem \ref{thm:all derivative},
with the role of the fundamental solution replaced by that of a positive
solution. The conditions \eqref{eq:2-9-2-4-1-1} guarantee that the
integration by parts formulas used for the fundamental solution in
the proof of Theorem \ref{thm:all derivative} also hold for the positive
solution. As in the proofs of Theorems \ref{thm:derivatives of TsallisEnt-1} and \ref{thm:Let--be},
for each $1\le i\le k$, the convergence of the integrals in \eqref{eq:2-9-2-5-1} and
\eqref{eq:2-9-2-6-1}, together with their continuity with respect to
$t$ on $(0,\infty)$, guarantees that these identities hold for every
$t>0$ under the continuity assumption in the theorem. 
\end{IEEEproof}
\subsection{Proofs of complete monotonicity and Euclidean consequences}

A real valued function $f(x)$ defined over an interval $I$ in the
real line is called a completely monotone function if its derivatives
satisfy $(-1)^{k}f^{(k)}(x)\ge0$ for each $k\ge0$ and all $x$ in
$I$. Under suitable geometry conditions, for the fundamental solution
of the heat equation, we show that the derivative $S_{2}^{\prime}$
is completely monotone with respect to $t$. Equivalently,
$(-1)^{k-1}S_{2}^{(k)}(u_t)\ge0$ for every integer $k\ge1$.
In the standard terminology, $S_{2}$ itself is not necessarily
completely monotone, since its first derivative is nonnegative.
\begin{IEEEproof}[Proof of Theorem~\ref{thm:fundamental-CM}]
This is a straightforward corollary of Theorem \ref{thm:all derivative}.
\end{IEEEproof}
For the positive solution of the heat equation, the corresponding
result is given as follows.
\begin{IEEEproof}[Proof of Theorem~\ref{thm:positive-CM}]
This is a straightforward corollary of Theorem \ref{thm:Let--be-1}.
\end{IEEEproof}
In the setting of Euclidean space and the classical heat equation,
the convergence and continuity properties of the integrals, as well
as the integration by parts formulas required in the computation of
the time derivatives of $S_{2}$, are satisfied. This follows from
the fact that a positive solution of the classical heat equation can
be represented as the convolution of its initial density with the
smooth Gaussian heat kernel. Thus, we obtain the following result.
\noindent\textit{Theorem~\ref{thm:For-any-positive-1} is an immediate Euclidean specialization of the $L^2$ spectral representation given below.}\par

The result in Theorem \ref{thm:For-any-positive-1} holds in Euclidean
spaces of arbitrary dimension; see \cite{wu2025onthecomp} for the
one-dimensional case.

\subsection{\label{sec:Spectrum-decomposition-for}Spectral proof and curvature-free representation}

\subsubsection{Representation by the spectral measure}

Let $P_t=e^{tL}$ be the symmetric Markov semigroup on $L^2(\mu)$,
and let $A=-L\ge0$ be its nonnegative self-adjoint generator. For
any $f\in L^2(\mu)$, the spectral theorem provides a positive
measure $\nu_f$ on $[0,\infty)$ such that
\[
\|P_tf\|_{L^2(\mu)}^2=\langle f,P_{2t}f\rangle
=\int_{[0,\infty)}e^{-2\lambda t}d\nu_f(\lambda).
\]
In particular, the corresponding $S_2(t)=1-\|P_tf\|_2^2$
satisfies, for every $k\ge1$ and $t>0$,
\[
S_2^{(k)}(t)=(-1)^{k-1}2^k
\int_{[0,\infty)}\lambda^k e^{-2\lambda t}d\nu_f(\lambda).
\]
The same representation holds for a heat kernel section
$p_t(\cdot,y)$ after choosing an arbitrary $s\in(0,t)$ and
setting $f=p_s(\cdot,y)$ (provided $f\in L^2(\mu)$), with
$t$ replaced by $t-s$. This spectral formulation does not require
a discrete spectrum.

\subsubsection{The discrete-spectrum specialization}

If $-L$ has a complete orthonormal eigenbasis $\{\varphi_k\}_{k\ge1}$
with eigenvalues $\lambda_k\ge0$, then under the usual $L^2$-convergence
assumption for a fixed heat-kernel section,
\begin{equation}
p_t(x,y)=\sum_{k\ge1}e^{-\lambda_kt}\varphi_k(x)\varphi_k(y).
\label{eq:decomp of p}
\end{equation}
The spectral theorem and Parseval's identity now give directly
\[
S_2(p_t(\cdot,y))=1-\sum_{k\ge1}e^{-2\lambda_kt}\varphi_k(y)^2,
\qquad
\frac{d}{dt}S_2(p_t(\cdot,y))
=2\sum_{k\ge1}\lambda_ke^{-2\lambda_kt}\varphi_k(y)^2.
\]
This is an $L^2$ calculation rather than a pointwise identity in $x$;
no discrete-spectrum assumption is required for the measure-theoretic
representation above. If $\mu$ is a probability measure and $-L$ has
a spectral gap $\lambda_{\mathrm{gap}}>0$ on the orthogonal complement
of constants, then the same spectral theorem yields
\[
\|P_tf-\mu(f)\|_2\le e^{-\lambda_{\mathrm{gap}}t}
\|f-\mu(f)\|_2.
\]
Finally, since a Laplace transform of a positive measure is log-convex,
$h_2(P_tf)=-\log\|P_tf\|_2^2$ is concave whenever the latter quantity
is finite. These are standard consequences of the spectral theorem
for symmetric Markov semigroups; see \cite{bakry2013analysis}.

\section{\label{sec:R=0000E9nyi-entropy-dissipation}Proofs of the R\'enyi entropy concavity theorems}

We compute the first and second time derivatives of Rényi entropy
$h_{\alpha}$ along the heat flow, where the order $\alpha>1$.

\subsection{Dissipation identities used in the concavity proof}

\begin{thm}
\label{thm:1st derivative of alpha}Under the assumptions of Theorem
\ref{thm:derivatives of TsallisEnt-1}, we have
\begin{equation}
\frac{d}{dt}h_{\alpha}(u)=\alpha\frac{\int_{M}|\nabla u(t)|^{2}u^{\alpha-2}(t)d\mu}{\int_{M}u^{\alpha}(t)d\mu},\label{eq:7-1}
\end{equation}
and
\begin{align}
\frac{d^{2}}{dt^{2}}h_{\alpha}(u) & =-\alpha\frac{2\int_{M}u^{\alpha-2}(Lu)^{2}d\mu+(\alpha-2)\int_{M}u^{\alpha-3}|\nabla u|^{2}Lud\mu}{\int_{M}u^{\alpha}(x)d\mu(x)}\nonumber \\
 & \quad-\alpha^{2}(1-\alpha)\left(\frac{\int_{M}|\nabla u(t)|^{2}u^{\alpha-2}(t)d\mu}{\int_{M}u^{\alpha}(x)d\mu(x)}\right)^{2}.\label{eq:7-2}
\end{align}
for every $t>0$.
\end{thm}
\begin{IEEEproof}
The formulas are direct consequences of the relation between
R\'enyi entropy and Tsallis entropy, a standard differentiation step
(see also \cite{wu2025onthecomp}). Indeed, writing
$Z_{\alpha}(t)=1+(1-\alpha)S_{\alpha}(u(t))=\int_Mu(t)^\alpha\,d\mu$,
we obtain
\begin{equation}
\frac{d}{dt}h_{\alpha}(u)=\frac{S_{\alpha}'(u)}{Z_{\alpha}(t)}.
\label{eq:7-3}
\end{equation}
and
\begin{equation}
\frac{d^{2}}{dt^{2}}h_{\alpha}(u)
=\frac{S_{\alpha}''(u)}{Z_{\alpha}(t)}
-(1-\alpha)\frac{(S_{\alpha}'(u))^{2}}{Z_{\alpha}(t)^{2}}.
\label{eq:7-4}
\end{equation}
Inserting the identities \eqref{eq:3-12-15} and \eqref{eq:2-9-2}
proves \eqref{eq:7-1} and \eqref{eq:7-2}. The exchange of differentiation
and integration is justified by Theorem~\ref{thm:derivatives of TsallisEnt-1}.
\end{IEEEproof}

We now prove the Rényi entropy concavity results stated in the Introduction.

\subsection{Proof for the fundamental solution}

\begin{IEEEproof}[Proof of Theorem~\ref{thm:Under-the-assumptions-1}]
Recall in Section \ref{sec:Preliminaries-for-the}, set $w=u^{\alpha/2}$,
$\alpha=\frac{2}{1+\delta}$,  then the first two derivative of Tsallis
entropy have the equivalent expressions \eqref{eq:3-6} and \eqref{eq:3-4}.
Note that the derivative expression \eqref{eq:7-4}, the sign of $\frac{d^{2}}{dt^{2}}h_{\alpha}(u)$
is always consistent with that of 
\begin{equation}
\frac{d^{2}}{dt^{2}}S_{\alpha}(u)\left(1+(1-\alpha)S_{\alpha}(u)\right)-(1-\alpha)\left(\frac{d}{dt}S_{\alpha}(u)\right)^{2}\label{eq:7-5}
\end{equation}
 Thus, we substitute \eqref{eq:3-6} and \eqref{eq:3-4} into the
derivative expression \eqref{eq:7-5}, we have
\begin{align*}
\eqref{eq:7-5} & =-4(1+\delta)\left(\int_{M}(Lw)^{2}d\mu+\delta\int_{M}Lw\frac{|\nabla w|^{2}}{w}d\mu\right)\int w^{2}d\mu-\frac{\delta-1}{1+\delta}\left(2(1+\delta)\int_{M}|\nabla w|^{2}d\mu\right)^{2}\\
 & =-4(1+\delta)\left(\left(\int_{M}(Lw)^{2}d\mu+\delta\int_{M}Lw\frac{|\nabla w|^{2}}{w}d\mu\right)\int w^{2}d\mu-(1-\delta)\left(\int_{M}|\nabla w|^{2}d\mu\right)^{2}\right).
\end{align*}
Since $\delta\in(-1,1]$, to prove that $\frac{d^{2}}{dt^{2}}h_{\alpha}(u)\le0$,
it suffices to verify the nonnegative of
\begin{equation}
\left(\int_{M}(Lw)^{2}d\mu+\delta\int_{M}Lw\frac{|\nabla w|^{2}}{w}d\mu\right)\int w^{2}d\mu-(1-\delta)\left(\int_{M}|\nabla w|^{2}d\mu\right)^{2}.\label{eq:7-6}
\end{equation}
Thanks to the integral estimate in Proposition \ref{prop:Under-the-assumption-1},
if $\delta\in[0,1]$ and $Ric(L)\ge0,$ then we substitute \eqref{eq:5-2}
into \eqref{eq:7-6},
\begin{equation}
\eqref{eq:7-6}\ge(1-\delta)\int_{M}(Lw)^{2}d\mu\int w^{2}d\mu-(1-\delta)\left(\int_{M}|\nabla w|^{2}d\mu\right)^{2}.\label{eq:9-3-1}
\end{equation}
By integration by parts formula \eqref{eq:4-1-1} in Proposition \ref{prop:For-.-we}
and Hölder's inequality,
\begin{align}
\left(\int_{M}|\nabla w|^{2}d\mu\right)^{2} & =\left(\int_{M}wLwd\mu\right)^{2}\nonumber \\
 & \le\int_{M}w^{2}d\mu\int_{M}(Lw)^{2}d\mu.\label{eq:9-3-2}
\end{align}
Substituting \eqref{eq:9-3-1} into \eqref{eq:9-3-2}, we have
\begin{align*}
\eqref{eq:7-6} & \ge(1-\delta)\int_{M}(Lw)^{2}d\mu\int w^{2}d\mu-(1-\delta)\left(\int_{M}|\nabla w|^{2}d\mu\right)^{2}\\
 & \ge(1-\delta)\int_{M}(Lw)^{2}d\mu\int w^{2}d\mu-(1-\delta)\int_{M}w^{2}d\mu\int_{M}(Lw)^{2}d\mu\\
 & =0.
\end{align*}
Therefore, when $\delta\in[0,1]$ and $Ric(L)\ge0,$ \eqref{eq:7-6}
is nonnegative. Finally, note that $\alpha=\frac{2}{1+\delta}$, thus
$\delta\in[0,1]$ corresponds to $\alpha\in[1,2]$, and we complete
the proof.
\end{IEEEproof}
\subsection{Proof for positive solutions and Euclidean specializations}

Under the assumptions of Theorem \ref{thm:Let--be}, we also obtain
the Rényi entropy dissipation formula for a positive solution of the heat
equation. 
\begin{IEEEproof}[Proof of Theorem~\ref{thm:Renyi-positive-concavity}]
The proof follows the same lines as that of Theorem \ref{thm:Under-the-assumptions-1},
with the fundamental solution replaced by a positive solution of the
heat equation.
\end{IEEEproof}
We also have the following counterpart in the Euclidean setting for
the classical heat equation.
\begin{IEEEproof}[Proof of Theorem~\ref{thm:For-any-positive}]
The proof is analogous to that of Theorem \ref{thm:Under-the-assumptions-1}.
In the Euclidean setting, we again set $w=u^{\alpha/2}$, $\alpha=\frac{2}{1+\delta}$
and rewrite the second derivative of Tsallis, and thus Rényi entropy
in terms of $w$. The integral estimates in Proposition \ref{prop:Under-the-assumption-1}
and the integration by parts formula in Proposition \ref{prop:For-.-we}
remain valid in Euclidean space; see Propositions 13, 14, and 16 in
\cite{sun2026concavity}. Therefore, the proof follows by applying
these results and Hölder’s inequality as in the proof of Theorem \ref{thm:Under-the-assumptions-1}.
\end{IEEEproof}
\begin{rem}
Our result in Theorem \ref{thm:For-any-positive} is new. In the one-dimensional
case, \cite{wu2025onthecomp} verified that when $\alpha\in(0,3]$,
the Rényi entropy $h_{\alpha}$ is concave in time along the heat
flow. The interval $[1,2]$ obtained in Theorem \ref{thm:For-any-positive}
is narrower than the interval established in the one-dimensional case.
Nevertheless, the interval $[1,2]$ is valid in Euclidean spaces of
arbitrary dimension. At present, we do not know whether the interval
$[1,2]$ in Theorem \ref{thm:Under-the-assumptions-1} or Theorem
\ref{thm:For-any-positive} is optimal. Extending this interval may
require the establishment of new functional inequalities and is left
for future investigation.
\end{rem}

\section{\label{sec:Concavity-of-R=0000E9nyi}Proofs of the R\'enyi entropy-power concavity theorems}

\subsection{Curvature-dimension reduction for the proof}

In $n$-dimensional Euclidean space $\mathbb{R}^{n}$, for $\alpha>0$
and $\beta>0$, the Rényi entropy power is defined by $e^{\frac{\beta}{n}h_{\alpha}(u)}$.
When $\alpha=1$ and $\beta=2$, it reduces to the classical Shannon
entropy power $e^{\frac{2}{n}h(u)}$. It is a classical result that
the Shannon entropy power is concave in time along the heat flow \cite{costa1985new},
see also \cite{li2024onthe} for the corresponding result on Riemannian
manifolds. However, as pointed out in Remark 2 of \cite{wu2025entropic},
when $\beta\ge2$, $e^{\frac{\beta}{n}h_{\alpha}(u)}$ can be concave
in time along the heat flow in Euclidean space only when $\alpha=1$.
Therefore, in order to obtain concavity results for Rényi entropy
powers that also cover the case $\alpha\ne1$, it is natural to consider
the range $0<\beta\le2$. Note that here we restrict our attention
to the (linear) heat flow. We refer to \cite{savare2014concavity}
for the concavity of Rényi entropy power along nonlinear heat flows,
and to \cite{li2020renyi} for the corresponding results on Riemannian
manifolds. 

\subsection{Proof for the fundamental solution}

Under the assumptions of Theorem \ref{thm:derivatives of TsallisEnt-1},
we define the Rényi entropy power by $e^{\frac{\beta}{m}h_{\alpha}(u)}$.
When $m=n$, it reduces to $e^{\frac{\beta}{n}h_{\alpha}(u)}$. We
then obtain the following result.
\begin{IEEEproof}[Proof of Theorem~\ref{thm:Under-the-assumptions-3}]
Denote $e^{\frac{\beta}{m}h_{\alpha}(u)}$ by $N_{\alpha,\beta}(u)$.
A direct computation yields
\begin{align*}
\frac{d}{dt}N_{\alpha,\beta}(u) & =\frac{\beta}{m}N_{\alpha,\beta}(u)\frac{d}{dt}h_{\alpha}(u),
\end{align*}
\begin{align*}
\frac{d^{2}}{dt^{2}}N_{\alpha,\beta}(u) & =N_{\alpha,\beta}(u)\left(\frac{\beta}{m}\frac{d}{dt}h_{\alpha}(u)\right)^{2}+\frac{\beta}{m}N_{\alpha,\beta}(u)\frac{d^{2}}{dt^{2}}h_{\alpha}(u)\\
 & =\frac{\beta}{m}N_{\alpha,\beta}(u)\left(\frac{\beta}{m}\left(\frac{d}{dt}h_{\alpha}(u)\right)^{2}+\frac{d^{2}}{dt^{2}}h_{\alpha}(u)\right).
\end{align*}
Substituting \eqref{eq:7-3} and \eqref{eq:7-4} in Theorem \ref{thm:1st derivative of alpha},
\begin{align*}
\frac{d^{2}}{dt^{2}}N_{\alpha,\beta}(u) & =\frac{\beta}{m}N_{\alpha,\beta}(u)\left(\frac{\beta}{m}\left(\frac{\frac{d}{dt}S_{\alpha}(u)}{1+(1-\alpha)S_{\alpha}(u)}\right)^{2}+\frac{\frac{d^{2}}{dt^{2}}S_{\alpha}(u)\left(1+(1-\alpha)S_{\alpha}(u)\right)-(1-\alpha)\left(\frac{d}{dt}S_{\alpha}(u)\right)^{2}}{\left(1+(1-\alpha)S_{\alpha}(u)\right)^{2}}\right)\\
 & =\frac{\beta}{m}N_{\alpha,\beta}(u)\frac{\left(\frac{\beta}{m}-1+\alpha\right)\left(\frac{d}{dt}S_{\alpha}(u)\right)^{2}+\frac{d^{2}}{dt^{2}}S_{\alpha}(u)\left(1+(1-\alpha)S_{\alpha}(u)\right)}{\left(1+(1-\alpha)S_{\alpha}(u)\right)^{2}}.
\end{align*}
Thus, to prove that $\frac{d^{2}}{dt^{2}}N_{\alpha,\beta}(u)\le0$,
it suffices to show that
\begin{equation}
\left(\frac{\beta}{m}-1+\alpha\right)\left(\frac{d}{dt}S_{\alpha}(u)\right)^{2}+\frac{d^{2}}{dt^{2}}S_{\alpha}(u)\left(1+(1-\alpha)S_{\alpha}(u)\right)\label{eq:9-1-1-1}
\end{equation}
is nonpositive. Recall from Section \ref{sec:Preliminaries-for-the}
that, by setting $w=u^{\alpha/2}$, $\alpha=\frac{2}{1+\delta}$,
and substituting \eqref{eq:3-6} and \eqref{eq:3-4} into \eqref{eq:9-1-1-1},
we can rewrite \eqref{eq:9-1-1-1} as
\[
4(1+\delta)\left(\left(\frac{\beta}{m}(1+\delta)+1-\delta\right)\left(\int_{M}|\nabla w|^{2}d\mu\right)^{2}-\left(\int_{M}(Lw)^{2}d\mu+\delta\int_{M}Lw\frac{|\nabla w|^{2}}{w}d\mu\right)\left(\int_{M}w^{2}d\mu\right)\right).
\]
Since $\delta\in(-1,1]$, to prove that $\frac{d^{2}}{dt^{2}}N_{\alpha,\beta}(u)\le0$,
it suffices to show that
\begin{equation}
\left(\frac{\beta}{m}(1+\delta)+1-\delta\right)\left(\int_{M}|\nabla w|^{2}d\mu\right)^{2}-\left(\int_{M}(Lw)^{2}d\mu+\delta\int_{M}Lw\frac{|\nabla w|^{2}}{w}d\mu\right)\left(\int_{M}w^{2}d\mu\right)\label{eq:9-1-1-2}
\end{equation}
is nonpositive. By the definition of $L$ and \eqref{eq:4-1-1} in
Proposition \ref{prop:For-.-we}, 
\[
\int_{M}w^{2}L\log wd\mu=\int_{M}wLwd\mu-\int_{M}|\nabla w|^{2}d\mu=2\int_{M}wLwd\mu.
\]
When $m=n$, the weight $\phi$ is constant by convention; any
term of the form $(m-n)^{-1}\langle\nabla\phi,\nabla w\rangle^{2}$
below is interpreted as zero. Then by the trace inequality,
\[
(L\log w)^{2}\le m|\nabla^{2}\log w|^{2}+\frac{m}{m-n}\langle\nabla\phi,\nabla\log w\rangle^{2},
\]
and it follows that
\begin{align}
4\left(\int_{M}wLwd\mu\right)^{2} & \le\int_{M}w^{2}d\mu\int_{M}w^{2}\left(L\log w\right)^{2}d\mu\nonumber \\
 & \le m\int_{M}w^{2}d\mu\int_{M}w^{2}|\nabla^{2}\log w|^{2}d\mu+\frac{m}{m-n}\int_{M}w^{2}d\mu\int_{M}\langle\nabla\phi,\nabla w\rangle^{2}d\mu.\label{eq:9-1-1-3}
\end{align}
 Combining inequality \eqref{eq:9-3-2} in the proof of Theorem \ref{thm:Under-the-assumptions-1}
with inequality \eqref{eq:9-1-1-3}, we obtain, for any $\theta\in[0,1]$,
\begin{align}
\left(\int_{M}|\nabla w|^{2}d\mu\right)^{2}
&\le\int_{M}w^{2}d\mu\left(\theta\int_{M}|Lw|^{2}d\mu\right.\nonumber\\
&\qquad\left.+\frac{(1-\theta)m}{4}\left(\int_{M}w^{2}|\nabla^{2}\log w|^{2}d\mu
+\frac{1}{m-n}\int_{M}\langle\nabla\phi,\nabla w\rangle^{2}d\mu\right)\right).
\label{eq:9-1-1-4}
\end{align}
When $0\le\delta\le1$ and
\[
\beta\ge\frac{4\delta-m(1-\delta)}{1+\delta}=2(2-\alpha)-m(\alpha-1),
\]
one can verify that
\begin{equation}
\theta=\frac{\beta(1+\delta)+m(1-\delta)-4\delta}{\beta(1+\delta)+m(1-\delta)}\label{eq:9-1-1-7}
\end{equation}
belongs to $[0,1]$, then substituting \eqref{eq:9-1-1-7} into \eqref{eq:9-1-1-4},
\begin{align}
\left(\frac{\beta(1+\delta)+m(1-\delta)}{m}\right)\left(\int_{M}|\nabla w|^{2}d\mu\right)^{2} & \le\int_{M}w^{2}d\mu\left(\frac{\beta(1+\delta)+m(1-\delta)-4\delta}{m}\int_{M}|Lw|^{2}d\mu\right.\nonumber \\
 & \quad\quad\left.+\delta\left(\int_{M}w^{2}|\nabla^{2}\log w|^{2}d\mu+\frac{1}{m-n}\int_{M}\langle\nabla\phi,\nabla w\rangle^{2}d\mu\right)\right).\label{eq:9-1-1-5}
\end{align}
For $0\le\delta\le1$, by \eqref{eq:2-9-5-1} in Proposition \ref{prop:Under-the-assumption-1},
\begin{align}
\eqref{eq:9-1-1-2} & =\left(\frac{\beta}{m}(1+\delta)+1-\delta\right)\left(\int_{M}|\nabla w|^{2}d\mu\right)^{2}\nonumber \\
 & \quad-\left(\delta\int_{M}w^{2}|\nabla^{2}\log w|^{2}d\mu+\delta\int_{M}\langle Ric(L)\nabla w,\nabla w\rangle d\mu+(1-\delta)\int_{M}|Lw|^{2}d\mu\right)\int_{M}w^{2}d\mu.\label{eq:9-1-1-6}
\end{align}
Substituting \eqref{eq:9-1-1-5} into \eqref{eq:9-1-1-6}, we derive
\begin{align*}
\eqref{eq:9-1-1-2} & \le\int_{M}w^{2}d\mu\left(\frac{\beta(1+\delta)+m(1-\delta)-4\delta}{m}\int_{M}|Lw|^{2}d\mu\right.\\
&\qquad\left.+\delta\left(\int_{M}w^{2}|\nabla^{2}\log w|^{2}d\mu+\frac{1}{m-n}\int_{M}\langle\nabla\phi,\nabla w\rangle^{2}d\mu\right)\right)\\
 & \quad-\left(\delta\int_{M}w^{2}|\nabla^{2}\log w|^{2}d\mu+\delta\int_{M}\langle Ric(L)\nabla w,\nabla w\rangle d\mu+(1-\delta)\int_{M}|Lw|^{2}d\mu\right)\int_{M}w^{2}d\mu\\
 & =\int_{M}w^{2}d\mu\left(\frac{\beta(1+\delta)-4\delta}{m}\int_{M}|Lw|^{2}d\mu-\delta\int_{M}\langle Ric_{m,n}(L)\nabla w,\nabla w\rangle d\mu\right).
\end{align*}
Thus, when $2(2-\alpha)-m(\alpha-1)\le\beta\le\frac{4\delta}{1+\delta}=2(2-\alpha)$,
$\beta>0$,  and 
\[
Ric_{m,n}(L)\ge0,
\]
we have
\[
\eqref{eq:9-1-1-2}\le0.
\]
It remains only to verify that if $\frac{d^{2}}{dt^{2}}N_{\alpha,\beta}(u)\le0$
for some $\beta>0$, then $\frac{d^{2}}{dt^{2}}N_{\alpha,\gamma}(u)\le0$
for every $0<\gamma<\beta$. Indeed, set $f(t)=N_{\alpha,\beta}(u)$,
$g(t)=N_{\alpha,\gamma}(u)$. A direct computation gives
\[
g'(t)=\left(f^{\frac{\gamma}{\beta}}(t)\right)'=\frac{\gamma}{\beta}f^{\frac{\gamma}{\beta}-1}(t)f'(t)
\]
and
\[
g''(t)=\frac{\gamma}{\beta}\left(\frac{\gamma}{\beta}-1\right)f^{\frac{\gamma}{\beta}-2}(t)(f'(t))^{2}+\frac{\gamma}{\beta}f^{\frac{\gamma}{\beta}-1}(t)f''(t).
\]
Since $\frac{\gamma}{\beta}<1$, the desired conclusion follows.
\end{IEEEproof}
\subsection{Proof for positive solutions and Euclidean specializations}

In the setting of Theorem \ref{thm:Under-the-assumptions-3}, we establish
the concavity of Rényi entropy power for the fundamental solution
of the heat equation. For a positive solution of the heat equation,
the corresponding concavity result is given as follows.
\begin{IEEEproof}[Proof of Theorem~\ref{thm:Renyi-positive-power}]
The proof follows the same lines as that of Theorem \ref{thm:Under-the-assumptions-3},
with the fundamental solution replaced by a positive solution of the
heat equation.
\end{IEEEproof}
We also establish the concavity of Rényi entropy power in the Euclidean
setting for the classical heat equation, which may be of particular
interest to the information theory community.
\begin{IEEEproof}[Proof of Theorem~\ref{thm:Renyi-euclidean-power}]
The proof is analogous to that of Theorem \ref{thm:Under-the-assumptions-3}.
In the Euclidean setting, we again set $w=u^{\alpha/2}$, $\alpha=\frac{2}{1+\delta}$,
and follow the same lines as in the proof of Theorem \ref{thm:Under-the-assumptions-3}
to rewrite the second derivative of Rényi entropy power in terms of
$w$. The integral estimates in Proposition \ref{prop:Under-the-assumption}
and the integration by parts formula in Proposition \ref{prop:For-.-we}
remain valid in Euclidean space; see Propositions 13 and 14 in \cite{sun2026concavity}.
Therefore, the proof follows by applying these results, together with
the trace inequality and Hölder’s inequality, as in the proof of Theorem
\ref{thm:Under-the-assumptions-3}.
\end{IEEEproof}

\appendix[Heat-kernel estimates used in the noncompact proofs]

Fix any point $y\in M$, let $p_{t}(x,y)$ denote the heat kernel
of the Witten Laplacian $L=\Delta-\nabla\phi\cdot\nabla$ with respect
to the weighted measure $d\mu=e^{-\phi}d\omega$. We denote by $cut(p)$
the cut locus of $M$ with respect to the point $p$.
\begin{lem}
\label{thm:gradient estimates-1}( \cite[Theorem 1.4]{li2016hamilton})
Let $(M,g)$ be an $n$-dimensional complete Riemannian manifold,
and $\phi\in C^{2}(M)$. Suppose that there exist constants $m\ge n$
and $K\ge0$ such that
\[
Ric_{m,n}(L)\ge-K.
\]
Then, for any fixed $T>0$, there exists a constant $C(K,m,n,T)>0$,
such that
\begin{equation}
|\nabla p_{t}(x,y)|\le C(K,m,n,T)\left(\frac{d(x,y)}{t}+\frac{1}{\sqrt{t}}\right)p_{t}(x,y)\label{eq:5-1-1}
\end{equation}
for all $x,y\in M$ and $t\in(0,T]$.
\end{lem}
The estimate is quoted from \cite[Theorem~1.4]{li2016hamilton};
we do not repeat its proof.
The following lemma gives the corresponding higher order gradient
estimates. 
\begin{lem}
\label{thm:gradient estimates} Let $N\ge2$. Let $(M,g)$ be an $n$-dimensional
complete Riemannian manifold. Assume that for every integer
$0\le j\le N$ there exists a constant $C^{(1)}_{j}>0$ such that
$\Vert\nabla^{j}Riem\Vert_{\infty}<C^{(1)}_{j}$. Let $\phi\in C^{N+1}(M)$
be such that $\nabla\phi\in C^{N}_{b}(M)$, and suppose that there
exist constants $m\ge n$ and $K\ge0$ such that
\[
Ric_{m,n}(L)\ge-K.
\]
Then, for any fixed $T>0$, there exists a constant $C_{N}(K,m,n,T)>0$,
such that
\begin{equation}
|\nabla^{N}p_{t}(x,y)|\le C_{N}(K,m,n,T)\left(\frac{d(x,y)}{t}+\frac{1}{\sqrt{t}}\right)^{N}p_{t}(x,y)\label{eq:5-1}
\end{equation}
for all $x,y\in M$ and $t\in(0,T]$.
\end{lem}
\begin{IEEEproof}
Put $F_t=\log p_t$ and $A_t=d(x,y)/t+t^{-1/2}$. By the
logarithmic heat-kernel derivative estimates of
\cite[Theorems~1.4 and~1.5]{li2016hamilton},
\begin{equation}
|\nabla^jF_t|\le C_jA_t^j,\qquad 1\le j\le N.
\label{eq:a2}
\end{equation}
The tensorial Fa\`a di Bruno formula gives
\begin{equation}
\nabla^Np_t=p_t\sum_{\pi\in\Pi_N}c_\pi\operatorname{Sym}
\left(\bigotimes_{B\in\pi}\nabla^{|B|}F_t\right).
\label{eq:a1}
\end{equation}
Every summand is bounded in norm by $C_\pi p_tA_t^N$ because
$\sum_{B\in\pi}|B|=N$. The finite sum yields \eqref{eq:5-1}.
This standard reduction from logarithmic derivative bounds to
heat-kernel derivative bounds is included only to identify precisely
the estimate used in the noncompact cutoff arguments.
\end{IEEEproof}
We also have the following useful lemma.
\begin{lem}
\label{lem:Let-.-Let} Let $(M,g)$ be an $n$-dimensional complete
Riemannian manifold. Let $\phi\in C^{2}(M)$, and suppose that there
exist constants $m\ge n$ and $K\ge0$ such that
\[
Ric_{m,n}(L)\ge-K.
\]
Then, for every nonnegative integer $q$, every fixed $y\in M$,
and every interval $(0,b]\subset(0,\infty),$ we have
\begin{equation}
\sup_{0<t\le b}\int_{M}d(x,y)^{q}p_{t}(x,y)d\mu(x)<\infty.\label{eq:a6}
\end{equation}
\end{lem}
\begin{IEEEproof}
We follow the argument used in the proof of Theorem 4.1 in \cite{li2012perelman}.
Fix $y\in M$, and write $r(x)=d(x,y)$. Let $\left(P_{t}\right)_{t\ge0}$
denote the heat semigroup generated by $L$. Thus, for every nonnegative
Borel function $f$,
\[
P_{t}f(x)=\mathbb{E}_{x}[f(X_{t})]=\int_{M}p_{t}(x,z)f(z)d\mu(z),
\]
where $\left(X_{t}\right)_{t\ge0}$ is the $L$-diffusion process
starting from $x$. Under the condition $Ric_{m,n}(L)\ge-K$, the
generalized Laplacian comparison theorem \cite{bakry2005volume,li2005liouville,qian1997estimates,qian1998comparison}
gives
\[
Lr\le(m-1)\sqrt{K}\coth(\sqrt{K}r)
\]
on $M\backslash(\{y\}\cup\text{Cut}(y))$. Using the inequality $a\coth a\le1+a,\forall a>0$,
we obtain
\begin{equation}
Lr\le(m-1)\left(\sqrt{K}+\frac{1}{r}\right)\label{eq:a5}
\end{equation}
on $M\backslash(\{y\}\cup cut(y))$. 

We first consider an integer $q\ge2$. Away from $y$ and its cut
locus, the distance function $r$ is smooth and satisfies $|\nabla r|=1$.
Consequently,

\[
Lr^{q}=qr^{q-1}Lr+q(q-1)r^{q-2}.
\]
By \eqref{eq:a5},
\begin{align*}
Lr^{q} & \le q(m-1)\sqrt{K}r^{q-1}+q(m-1)r^{q-2}+q(q-1)r^{q-2}\\
 & =q(m-1)\sqrt{K}r^{q-1}+q(m+q-2)r^{q-2}.
\end{align*}
Since $q\ge2$,
\[
r^{q-1}\le1+r^{q},\quad r^{q-2}\le1+r^{q},
\]
for every $r\ge0$. It follows that there exist constants
\[
C_{1}=C_{1}(q,m,K)>0,\quad C_{2}=C_{2}(q,m,K)>0,
\]
such that
\begin{equation}
Lr^{q}\le C_{1}r^{q}+C_{2}\label{eq:a7}
\end{equation}
on $M\backslash(\{y\}\cup\text{Cut}(y))$. 

Set $r(X_{t})=d(X_{t},y)$. As in the proof of Theorem 4.1 in \cite{li2012perelman},
there exists an one-dimensional Brownian motion $\beta_{t}$ such
that \cite{cranston1991gradient,kendall1986nonnegative,li2005liouville,wang1997logarithmic}
\[
dr(X_{t})\le d\beta_{t}+Lr(X_{t})dt.
\]
Applying Itô's formula to $r^{q}(X_{t})$, with the usual localization
understood, and using \eqref{eq:a7}, we obtain
\begin{align*}
dr^{q}(X_{t}) & \le qr^{q-1}(X_{t})d\beta_{t}+Lr^{q}(X_{t})dt\\
 & \le qr^{q-1}(X_{t})d\beta_{t}+C_{1}r^{q}(X_{t})dt+C_{2}dt.
\end{align*}
Taking expectation with respect to the law of $X_{t}$,
\[
\frac{d}{dt}\mathbb{E}_{x}[r^{q}(X_{t})]\le C_{1}\mathbb{E}_{x}[r^{q}(X_{t})]+C_{2}.
\]
Equivalently,
\[
\mathbb{E}_{x}[r^{q}(X_{t})]\le r(x)^{q}+C_{2}t+C_{1}\int^{t}_{0}\mathbb{E}_{x}[r^{q}(X_{s})]ds.
\]
By Grönwall inequality,
\[
\mathbb{E}_{x}[r^{q}(X_{t})]\le e^{C_{1}t}\left(r(x)^{q}+C_{2}t\right).
\]
In terms of the heat semigroup, this can be written as
\begin{equation}
P_{t}\left[d(\cdot,y)^{q}\right](x)\le e^{C_{1}t}\left(d(x,y)^{q}+C_{2}t\right).\label{eq:2-9-9}
\end{equation}
Taking $x=y$ in \eqref{eq:2-9-9}, we obtain
\[
P_{t}\left[d(\cdot,y)^{q}\right](y)\le C_{2}te^{C_{1}t}.
\]
Since $L$ is symmetric with respect to $\mu$, the heat kernel is
symmetric, i.e.,
\[
p_{t}(x,y)=p_{t}(y,x).
\]
Hence,
\begin{align*}
P_{t}\left[d(\cdot,y)^{q}\right](y) & =\int_{M}d(x,y)^{q}p_{t}(y,x)d\mu(x)\\
 & =\int_{M}d(x,y)^{q}p_{t}(x,y)d\mu(x).
\end{align*}
Therefore,
\[
\sup_{0<t\le b}\int_{M}d(x,y)^{q}p_{t}(x,y)d\mu(x)\le C_{2}be^{C_{1}b}<\infty
\]
for every integer $q\ge2$.

For $q=1$, the conservativeness of the heat semigroup and the Cauchy-Schwarz
inequality give
\begin{align*}
P_{t}\left[d(\cdot,y)\right](y) & \le\left(P_{t}\left[d(\cdot,y)^{2}\right](y)\right)^{\frac{1}{2}}\left(P_{t}1(y)\right)^{\frac{1}{2}}\\
 & =\left(P_{t}\left[d(\cdot,y)^{2}\right](y)\right)^{\frac{1}{2}},
\end{align*}
which is uniformly bounded for $0<t\le b$. Finally, for $q=0$,
\[
P_{t}1(y)=\int_{M}p_{t}(x,y)d\mu(x)=1.
\]
Thus, the conclusion follows.

\end{IEEEproof}
%

% === End of references ===


\begin{thebibliography}{10}

\bibitem{bakry2006diffusions}
Dominique Bakry and Michel {\'E}mery.
\newblock Diffusions hypercontractives.
\newblock In {\em S{\'e}minaire de Probabilit{\'e}s XIX 1983/84: Proceedings},
  pages 177--206. Springer, 1985.

\bibitem{bakry2013analysis}
Dominique Bakry, Ivan Gentil, and Michel Ledoux.
\newblock {\em Analysis and Geometry of {Markov} Diffusion Operators}, volume
  348 of {\em Grundlehren der mathematischen Wissenschaften}.
\newblock Springer, 2014.

\bibitem{bakry2005volume}
Dominique Bakry and Zhongmin Qian.
\newblock Volume comparison theorems without jacobi fields.
\newblock In {\em Conference on Potential Theory}, pages p--115, 2005.

\bibitem{Cheng2015Higher}
Fan Cheng and Yanlin Geng.
\newblock Higher order derivatives in {Costa}'s entropy power inequality.
\newblock {\em IEEE Transactions on Information Theory}, 61(11):5892--5905,
  2015.

\bibitem{costa1985new}
M.~Costa.
\newblock A new entropy power inequality.
\newblock {\em IEEE Transactions on Information Theory}, 31(6):751--760, 1985.

\bibitem{cranston1991gradient}
Michael Cranston.
\newblock Gradient estimates on manifolds using coupling.
\newblock {\em Journal of functional analysis}, 99(1):110--124, 1991.

\bibitem{dembo1989simple}
Amir Dembo.
\newblock Simple proof of the concavity of the entropy power with respect to
  added {Gaussian} noise.
\newblock {\em IEEE Transactions on Information Theory}, 35(4):887--888, 1989.

\bibitem{fisher1925theory}
Ronald~Aylmer Fisher.
\newblock Theory of statistical estimation.
\newblock {\em Mathematical Proceedings of the Cambridge Philosophical
  Society}, 22:700--725, 1925.

\bibitem{gross1975logarithmic}
Leonard Gross.
\newblock Logarithmic {Sobolev} inequalities.
\newblock {\em American Journal of Mathematics}, 97(4):1061--1083, 1975.

\bibitem{guo2021generalization}
Laigang Guo, Chuang-Ming Yuan, and Xiao-Shan Gao.
\newblock A generalization of the concavity of {R{\'e}nyi} entropy power.
\newblock {\em Entropy}, 23(12):1593, 2021.

\bibitem{hung2022generalization}
Li-Chang Hung.
\newblock Generalization of completely monotone conjecture for {Tsallis}
  entropy, 2022.
\newblock arXiv:2212.09269.

\bibitem{kendall1986nonnegative}
Wilfrid~S Kendall.
\newblock Nonnegative ricci curvature and the brownian coupling property.
\newblock {\em Stochastics: An International Journal of Probability and
  Stochastic Processes}, 19(1-2):111--129, 1986.

\bibitem{ledoux1994semigroup}
Michel Ledoux.
\newblock Semigroup proofs of the isoperimetric inequality in {Euclidean} and
  {Gauss} space.
\newblock {\em Bulletin des Sciences Math{\'e}matiques}, 118(6):485--510, 1994.

\bibitem{ledoux2022differentials}
Michel Ledoux.
\newblock Differentials of entropy and {Fisher} information along heat flow: A
  brief review of some conjectures, 2022.

\bibitem{ledoux2021log}
Michel Ledoux, Chandra Nair, and Yan~Nan Wang.
\newblock Log-convexity of {Fisher} information along heat flow, 2021.

\bibitem{li2020renyi}
S.~Li and X.-D. Li.
\newblock On the {R{\'e}nyi} entropy power and the
  {Gagliardo-Nirenberg-Sobolev} inequality on {Riemannian} manifolds.
\newblock {\em arXiv preprint arXiv:2001.11184}, 2020.

\bibitem{li2024onthe}
Songzi Li and Xiang-dong Li.
\newblock On the {Shannon} entropy power on {Riemannian} manifolds and {Ricci}
  flow.
\newblock {\em Tohoku Mathematical Journal}, 76(4):577--608, 2024.

\bibitem{li2005liouville}
Xiang-Dong Li.
\newblock Liouville theorems for symmetric diffusion operators on complete
  riemannian manifolds.
\newblock {\em Journal de math{\'e}matiques pures et appliqu{\'e}es},
  84(10):1295--1361, 2005.

\bibitem{li2012perelman}
Xiang-Dong Li.
\newblock {Perelman}’s entropy formula for the {Witten} laplacian on
  {Riemannian} manifolds via {Bakry--Emery} {Ricci} curvature.
\newblock {\em Mathematische Annalen}, 353(2):403--437, 2012.

\bibitem{li2016hamilton}
Xiang-Dong Li.
\newblock {Hamilton}’s {Harnack} inequality and the {$W$}-entropy formula on
  complete {Riemannian} manifolds.
\newblock {\em Stochastic Processes and their Applications}, 126(4):1264--1283,
  2016.

\bibitem{lott2002some}
John Lott.
\newblock Some geometric properties of the bakry-{\'e}mery-ricci tensor.
\newblock {\em Comment. Math. Helv.}, 78:865--883, 2003.

\bibitem{McKean1966}
Henry~P. McKean.
\newblock Speed of approach to equilibrium for {Kac}'s caricature of a
  {Maxwellian} gas.
\newblock {\em Archive for Rational Mechanics and Analysis}, 21(5):343--367,
  1966.

\bibitem{qian1997estimates}
Zhongmin Qian.
\newblock Estimates for weighted volumes and applications.
\newblock {\em Quarterly Journal of Mathematics}, 48(2):235--242, 1997.

\bibitem{qian1998comparison}
Zhongmin Qian.
\newblock A comparison theorem for an elliptic operator.
\newblock {\em Potential Analysis}, 8(2):137--142, 1998.

\bibitem{Renyi1961OnMO}
Alfr{\'e}d R{\'e}nyi.
\newblock On measures of entropy and information.
\newblock In {\em Proceedings of the Fourth Berkeley Symposium on Mathematical
  Statistics and Probability}, volume~1, pages 547--561. University of
  California Press, 1961.

\bibitem{savare2014concavity}
G.~Savar{\'e} and G.~Toscani.
\newblock The concavity of {R{\'e}nyi} entropy power.
\newblock {\em IEEE Transactions on Information Theory}, 60(5):2687--2693,
  2014.

\bibitem{Shannon48}
C.~E. Shannon.
\newblock {A Mathematical Theory of Communication}.
\newblock {\em Bell System Technical Journal}, 27:379--423, 1948.

\bibitem{StamEPI59}
A.~J. Stam.
\newblock Some inequalities satisfied by the quantities of information of
  {Fisher} and {Shannon}.
\newblock {\em Information and Control}, 2:101--112, 1959.

\bibitem{sun2026concavity}
Lukang Sun.
\newblock Concavity of {Tsallis} entropy and {Tsallis} entropy power along heat
  flow, 2026.
\newblock arXiv:2601.12944, revised version.

\bibitem{villani2000short}
C{\'e}dric Villani.
\newblock A short proof of the ``concavity of entropy power''.
\newblock {\em IEEE Transactions on Information Theory}, 46(4):1695--1696,
  2000.

\bibitem{wang1997logarithmic}
Feng-Yu Wang.
\newblock Logarithmic sobolev inequalities on noncompact riemannian manifolds.
\newblock {\em Probability theory and related fields}, 109(3):417--424, 1997.

\bibitem{wu2025entropic}
H.~Wu and L.~Yu.
\newblock Entropic isoperimetric and {Cramér--Rao} inequalities for
  {Rényi--Fisher} information.
\newblock {\em IEEE Transactions on Information theory}, 72(2):901--920, 2026.

\bibitem{wu2025onthecomp}
H.~Wu, L.~Yu, and L.~Guo.
\newblock On the complete monotonicity of {R{\'e}nyi} entropy.
\newblock {\em IEEE Transactions on Information Theory}, 71(7):4895--4914,
  2025.

\end{thebibliography}
\end{document}